\documentclass[a4paper,12pt,french,leqno]{smfart}

\usepackage{hyperref}
\usepackage{amsmath}
\usepackage{amssymb}
\usepackage{amsthm}
\usepackage{amsopn}
\usepackage{amscd}
\usepackage{fullpage}
\usepackage[utf8]{inputenc}
\usepackage{graphics, xcolor}
\usepackage{textcomp} 
\usepackage[all]{xy}
\usepackage{lscape}
\usepackage{url}
\usepackage{verbatim}
\usepackage{mathrsfs}
  \usepackage[francais]{babel}

\usepackage{bm}

\def\commentaire#1{\ifx\cachercommentaires\undefined \textcolor{red}{#1}\else \fi } 

\newcommand{\scrA}{\mathscr{A}}

\newcommand{\scrF}{\mathscr{F}}
\newcommand{\scrG}{\mathscr{G}}

\newcommand{\scrK}{\mathscr{K}}

\newcommand{\frU}{\mathfrak{U}}

\newcommand{\frZ}{\mathfrak{Z}}

\newcommand{\frb}{\mathfrak{b}}

\newcommand{\frg}{\mathfrak{g}}
\newcommand{\frh}{\mathfrak{h}}

\newcommand{\frk}{\mathfrak{k}}
\newcommand{\frl}{\mathfrak{l}}
\newcommand{\frmm}{\mathfrak{m}}
\newcommand{\frn}{\mathfrak{n}}

\newcommand{\frp}{\mathfrak{p}}
\newcommand{\frqqq}{\mathfrak{q}}

\newcommand{\frt}{\mathfrak{t}}
\newcommand{\fru}{\mathfrak{u}}
\newcommand{\frv}{\mathfrak{v}}

\newcommand{\frz}{\mathfrak{z}}

\newcommand{\bbC}{\mathbb{C}}

\newcommand{\bbR}{\mathbb{R}}

\newcommand{\bbZ}{\mathbb{Z}}

\newcommand{\caD}{\mathcal{D}}
\newcommand{\caE}{\mathcal{E}}

\newcommand{\caM}{\mathcal{M}}
\newcommand{\caN}{\mathcal{N}}

\newcommand{\caR}{\mathcal{R}}

\newcommand{\caT}{\mathcal{T}}

\def\U{\mathbf{U}}

\newcommand{\GL}{\mathbf{GL}}

\newcommand{\Ad}{\mathrm{Ad}}

\newcommand{\End}{\mathrm{End}}

\newcommand{\Ind}{\mathrm{Ind}}

\renewcommand{\Ind}{\mathrm{Ind}}

\newcommand{\SL}{\mathbf{SL}}

\newcommand{\SO}{\mathbf{SO}}
\newcommand{\Sp}{\mathbf{Sp}}

\newcommand{\bil}[2]{\langle  #1,#2 \rangle }

\newcommand{\sgn}{\mathbf{sgn}}
\newcommand{\Triv}{\mathbf{Triv}}

\newcommand{\Std}{\mathbf{Std}}

\theoremstyle{plain}
\newtheorem{thm}{Théorème}[section]
\newtheorem{lemme}[thm]{Lemme}

\newtheorem{prop}[thm]{Proposition}

\theoremstyle{definition}
 \newtheorem{defi}[thm]{Définition}

\newtheorem{rmq}[thm]{Remarque}

\def \dem {\noindent \underline{\sl Démonstration}. }

\begin{document}

\numberwithin{equation}{section}

\title{Sur les paquets d'Arthur des groupes unitaires et quelques conséquences pour les groupes classiques}
  \author{Colette Moeglin}
 \address{CNRS, Institut Mathématique de Jussieu } 
 \email{colette.moeglin@imj-prg.fr}

  \author{David Renard  }
 \address{Centre de Mathématiques
  Laurent Schwartz,  Ecole Polytechnique} 
\email{david.renard@polytechnique.edu}

\date{\today}


\begin{abstract} Nous donnons une construction explicite des paquets d'Arthur des groupes unitaires réels 
par induction cohomologique et induction parabolique
 et en suivant une idée communiquée par P. Trapa, nous établissons la propriété de multiplicité un de ceux-ci. Nous montrons en particulier 
des résultats d'irréductibilité de certaines induites paraboliques pour les groupes unitaires, ce qui nous permet de compléter
les démonstrations d'énoncés analogues  annoncés  dans nos travaux sur les paquets d'Arthur des groupes classiques. 
\medskip
 
\noindent {\bf  Abstract.}\, ---\, We give an explicit construction of Arthur packets for real unitary groups by cohomological and parabolic induction
and following an idea communicated to us by P. Trapa,  we show that they satisfy the  multiplicity one property. In particular, we show the irreducibility 
of some parabolically induced representations
for unitary groups, and use this to give the proof of analogous statements made in our work on Arthur packets of classical
groups.
\end{abstract}

\bigskip
\maketitle

\section{Introduction} Le premier  objet  de cet article est de déterminer explicitement les paquets d'Arthur des groupes unitaires réels, et d'établir 
un résultat de multiplicité un pour ceux-ci. Dans \cite{MR3} et les articles afférents \cite{MR4}, \cite{MR5}, des résultats analogues
ont été établis pour les groupes classiques  ({\sl ie.} spéciaux orthogonaux et symplectiques) réels. Nous complétons  aussi nos résultats sur les groupes
classiques en donnant les démonstrations d'énoncés de réduction aux paramètres de bonne parité et d'irréductibilité d'induites paraboliques
dans cette réduction  annoncés dans \cite{MR3}. Cette  démonstration d'irréductibilité d'induites pour les groupes classiques utilise le résultat
analogue  pour les groupes unitaires démontré  dans cet article, ce qui explique qu'elle apparaisse seulement  ici.  

Les paquets d'Arthur des groupes classiques et unitaires sont déterminés par leurs propriétés, plus précisément
par certaines identités endoscopiques (\cite{Art13}, \cite{Mok}). Le but est de donner une construction la plus explicite possible
des représentations dans ces paquets. Nous renvoyons à \cite{MR3}, \S 2 pour une discussion générale sur les paquets d'Arthur
et  \S 4 pour les énoncés de nos résultats pour les groupes classiques. Rappelons simplement  ici quelques éléments.
Soit $\mathbf G$ un groupe algébrique  connexe réductif défini sur $\bbR$ et notons $G$ le groupe de ses points réels. 
Soit $\psi_G: \, W_\bbR \times \SL_2(\bbC) \rightarrow {}^LG$ un paramètre d'Arthur. Notons $S_{\psi_G}$ le centralisateur de l'image de 
$\psi_{G}$ dans $\widehat G$,  $(S_{\psi_G})_0$   sa composante neutre, et posons $A(\psi_G)=S_{\psi_G}/(S_{\psi_G})_0$. Supposons pour simplifier que 
ces groupes soient abéliens (c'est une hypothèse qui porte sur $\mathbf G$ et qui est vérifiée si $\mathbf G$ est un groupe classique ou unitaire).  
D'autre part, supposons que $\mathbf G$ soit quasi-déployé, ou bien forme intérieure pure d'un groupe quasi-déployé.
Les conjectures d'Arthur (\cite{Art84}, \cite{Art89}) dans ce cadre reviennent   alors à affirmer  l'existence d'une certaine combinaison linéaire à coefficients complexes 
de représentations irréductibles unitaires de $G\times A(\psi_G)$,  que nous notons $\pi^A(\psi_G)$, et qui doit vérifier certaines propriétés, notamment 
des identités de transfert endoscopique. Pour les groupes classiques quasi-déployés, ceci est établi dans \cite{Art13}, où $\pi^A(\psi_G)$
est caractérisée par les  identités de transfert endoscopique attachées aux données endoscopiques elliptiques de $\mathbf G$ et par une identité
 de transfert endoscopique vers un groupe général linaire tordu. Il est démontré de plus 
que les coefficients complexes dans $ \pi^A(\psi_G)$ sont en fait des entiers positifs. Ainsi on peut voir  $ \pi^A(\psi_G)$ comme une représentation 
unitaire de longueur finie de $G\times A(\psi_G)$. Décomposons cette représentation selon les représentations irréductibles unitaires de $G$
en écrivant 
\[ \pi^A(\psi_G)=\bigoplus_{\pi\in \Pi(\psi_G)} \pi \boxtimes \rho_\pi .\]
Ici, $\Pi(\psi_G)$ est donc un ensemble fini de  représentations irréductibles unitaires de $G$, le paquet d'Arthur attaché à $\psi_G$, et pour tout $\pi \in \Pi(\psi_G)$, 
$\rho_\pi$ est une représentation unitaire de dimension finie de $A(\psi_G)$ (une somme directe de caractères car $A(\psi_G)$ est abélien).
La dimension de $\rho_\pi$ est la multiplicité de $\pi$ dans le  paquet $\Pi(\psi_G)$ (rappelons que les paquets d'Arthur ne sont pas disjoints).
Ces résultats sont aussi démontrés pour les  groupes classiques non quasi-déployés, c'est-à-dire les groupes spéciaux orthogonaux réels $\SO(p,q)$, 
dans \cite{MR4} et pour les groupes  groupes unitaires  dans \cite{Mok} et \cite{KMSW}.

Le travail entrepris dans \cite{MR3} auquel nous renvoyons pour plus de détails, est de donner une construction explicite de  $\pi^A(\psi_G)$ pour les groupes classiques.
Un problème important qui nous occupe  aussi est d'établir  que les multiplicités sont $1$, c'est-à-dire que $\rho_\pi$ est un caractère de $A(\psi_G)$.
Considérons la composition de $\psi_G$ avec la représentation standard $\Std_{\mathbf G}$ du $L$-groupe de $\mathbf G$. On obtient un 
paramètre 
\[\psi={\Std_\mathbf G}\circ \psi_G : W_\bbR \times \SL_2(\bbC) \longrightarrow \GL_N(\bbC)\]
 que l'on voit comme une représentation complètement réductible de dimension $N$ 
de $W_\bbR\times \SL_2(\bbC)$. On écrit une décomposition de $\psi$ de la forme 
\[\psi=\psi_{mp}\oplus\psi_{bp}= \psi_{mp}\oplus\psi_{bp,disc}\oplus  \psi_{bp,u}\]
où  $\psi_{mp}$ est la partie de mauvaise parité du paramètre, $\psi_{bp,disc}$ la partie de bonne parité discrète, et 
$\psi_{bp,u}$ la partie de bonne parité unipotente.  Les constructions  de \cite{MR3} se font en quatre étapes.
La première étape est le cas $\psi=\psi_{bp,u}$ des paramètres unipotents et de bonne parité.
Dans ce cas, $\pi^A(\psi_G)$ est déterminée dans \cite{pourhowe} par des correspondance de Howe
et la propriété de multiplicité un de tels  paquets y est établie. Les représentations de $G$  dans ces paquets sont faiblement unipotentes au sens de \cite{KnVo}, 
chapter XII.  La deuxième étape est le cas où $\psi= \psi_{bp,disc}\oplus  \psi_{bp,u}$, lorsque le paramètre $ \psi_{bp,disc}$ possède certaines propriétés de régularité.
Les représentations dans le paquet $\Pi(\psi_G)$ sont alors obtenues par l'induction cohomologique de Vogan-Zuckerman
à partir des représentations faiblement unipotentes dans le paquet  $\Pi(\psi_{G_u})$, où 
$\psi_{G_u}$ est un paramètre d'Arthur unipotent pour un groupe $\mathbf G_u$ de même type que $\mathbf G$ et de rang plus petit,
qui après composition avec la représentation standard $\Std_{\mathbf G_u}$ donne  $\psi_{bp,u}$,  et de caractères de groupes unitaires
associés à $\psi_{u,disc}$.
Sous l'hypothèse de régularité mentionnée, les inductions cohomologiques se font dans le \og good range \fg, et en particulier sont irréductibles, 
les paramètres de Langlands des induites se déduisent facilement de ceux des induisantes, et la propriété de multiplicité un des paquets est conservée.
Les résultats sont  démontrés dans \cite{MR3} et utilisent de manière cruciale les constructions d'Adams et Johnson \cite{AdJo87}, \cite{Joh1} qui sont reliées à 
celles d'Arthur dans  \cite{AMR}.
La troisième étape établie dans \cite{MR5}  consiste à s'affranchir de l'hypothèse de régularité de $ \psi_{bp,disc}$,
 et l'on utilise pour cela les propriétés des foncteurs de translation 
(\cite{KnVo}, chapter VII). Les représentations dans   $\Pi(\psi_G)$  sont encore obtenues par induction cohomologique comme dans le cas régulier, mais 
celle-ci a maintenant lieu dans le  \og weakly fair range \fg, où les résultats généraux sont moins fort. En particulier, l'irréductibilité n'est plus préservée, il peut y avoir
 des annulations, les paramètres de Langlands des induites deviennent très délicats à calculer car  il n'y a pas de formule générale.
En conséquence, on perd dans cette étape la conservation de la propriété de multiplicité un (mais nous conjecturons que celle-ci reste vraie, il faudrait l'établir 
en utilisant d'autres outils).
 Enfin, la quatrième étape consiste à passer des paquets de bonne
parité aux paquets généraux. Les énoncés sont simples, cela se fait par une induction parabolique 
à partir des représentations du paquet $\Pi(\psi_{G_{bp}})$, où  
$\psi_{G_{bp}}$ est un paramètre d'Arthur de bonne parité pour un groupe $\mathbf G_{bp}$ de même type que $\mathbf G$ et de rang plus petit,
qui après composition avec la représentation standard $\Std_{\mathbf G_{bp}}$ donne  $\psi_{bp}$ et d'une représentation d'un groupe linéaire
attachée à $\psi_{mp}$. Ces résultats de réduction à la bonne parité ont  été annoncés dans \cite{MR3} sans démonstrations, et
nous donnons  celles-ci ici (Proposition \ref{reducPaqArt} et Théorème \ref{redCla}).

L'objet principal de cet article est donc d'établir des résultats analogues, mais pour les groupes unitaires. La stratégie suit les mêmes étapes, 
mais les résultats sont plus simples. Soit $\mathbf G$ un groupe unitaire de rang $N$  défini sur $\bbR$,
disons $G=\U(p,q)$, $p+q=N$, 
 et $\psi_G: W_\bbR\times \SL_2(\bbC)\rightarrow {}^LG$
un paramètre d'Arthur pour $G$. Maintenant, on remplace la composition avec la représentation standard par la restriction du paramètre à 
$\bbC^\times \hookrightarrow W_\bbR$ (changement de base). On obtient donc un paramètre 
\[\psi=\psi_{G\vert \bbC^\times}: \bbC^\times\times \SL_2(\bbC) =W_\bbC\times \SL_2(\bbC) \longrightarrow \GL_N(\bbC),\]
que l'on voit comme un paramètre pour $\GL_N(\bbC)$. Là encore, on décompose $\psi$ en somme de représentations irréductibles de $\bbC^\times$
(donc de dimension $1$) et l'on sépare les composantes de bonne et de mauvaise parité  (voir section \ref{decAparU}) :
\[\psi=\psi_{mp}\oplus \psi_{bp}.\]
La première simplification par rapport au cas des groupes classiques et qu'il n'y a pas à considérer de partie unipotente.
La première étape consiste donc à établir les résultats dans le cas $\psi=\psi_{bp}$, en faisant d'abord   là aussi une hypothèse de régularité du paramètre.
Ecrivons 
\[\psi=\psi_{bp}=\bigoplus_{i=1}^\ell ( \chi_{t_i}\boxtimes R[a_i])\]
où $R[a_i]$ est la représentation irréductible algébrique de dimension $a_i$ de $\SL_2(\bbC)$ et $\chi_{t_i}$, $t_i\in \bbZ$,   est le caractère de $\bbC^\times$
défini par  
$z\mapsto \left(\frac{z}{\bar z}\right)^{\frac{t_i}{2}}$. La condition de bonne parité est que pour tout $i=1,\ldots,\ell$, $t_i+a_i-N$ est pair.
On  suppose les $t_i$ rangés dans l'ordre décroissant, ce qui est loisible. La condition de régularité est alors  que pour tout $i=1,\ldots,\ell-1$, 
\begin{equation}\label{regintro} t_i-(a_i-1)>t_{i+1}+(a_{i+1}-1).\end{equation}
Remarquons que $N=\sum_{i=1}^\ell a_i$. On note $\caD(\psi)$ l'ensemble des familles $\underline d =(p_i,q_i)_{i=1,\ldots,\ell}$ de couples d'entiers tels que 
$\sum_{i=1}^\ell p_i=p$ et $\sum_{i=1}^\ell q_i=q$.
La représentation $\pi^A(\psi)$ est construite par induction cohomologique. Notons $\frg$ l'algèbre de Lie complexifiée de $G=\U(p,q)$, 
$K$ un sous-groupe compact maximal de $G$ associé à une involution de Cartan $\theta$, et $\frk$ la complexification de l'algèbre de Lie de $\frk$.
A un élément $\underline d$ de $\caD(\psi)$, on associe de manière explicite une sous-algèbre parabolique 
  $\theta$-stable $\frqqq_{\underline d}=\frl_{\underline d}
\oplus \frv_{\underline d}$ de $\frg$. On pose $L_{\underline d}=\mathrm{Norm}_G(\frqqq_{\underline d})$. C'est un $c$-Levi de $G$, au sens de Shelstad
\cite{She15}, isomorphe au produit $\times_{i=1}^\ell \U(p_i,q_i)$, et la complexifiée de l'algèbre de Lie de ce groupe est $\frl_{\underline d}$. 
On introduit un caractère $\Lambda_{\underline d}$ de ce groupe, explicitement déterminé par les $(t_i,a_i)$ (voir Equation (\ref{lambdads})), et l'on pose 
\begin{equation}\label{Adintro}
\scrA_{\underline d}(\psi)=\left(\caR_{\frqqq_{\underline d}, L_{\underline d}\cap K}^{\frg,K}\right)^{\dim(\frv_{\underline d}\cap \frk)} (\Lambda_{\underline d})
\end{equation}
où $\left(\caR_{\frqqq_{\underline d}, L_{\underline d}\cap K}^{\frg,K}\right)^{k}$ est le foncteur d'induction cohomologique de Vogan-Zuckerman en degré $k$ ({\sl cf.} \cite{KnVo} chapter V).
La condition  de régularité (\ref{regintro}) assure que cette induction cohomologique est dans le  \og good  range \fg.
De ceci, il découle que $\left(\caR_{\frqqq_{\underline d}, L_{\underline d}\cap K}^{\frg,K}\right)^{k} (\Lambda_{\underline d})=0$,
 si  $k\neq \dim(\frv_{\underline d}\cap \frk)$ , et que   $\scrA_{\underline d}(\psi)$ est un module unitaire
et irréductible. 
 De plus, les $\scrA_{\underline d}(\psi)$ lorsque $\underline d$ parcourt $\caD(\psi)$ sont distincts.

Définissons maintenant  pour tout $\underline d\in \caD(\psi)$ un caractère  $\epsilon_{\underline{d}}$ du groupe $A(\psi)$.
On définit d'abord $\epsilon_{\underline{d}}$ comme une application de $[1,\ell]$ dans $\pm 1$. Pour cela on pose pour tout entier $ i\in [1,\ell]$,
$a_{<i}=\sum_{j<i} a_j$ et 
\begin{equation}\label{edintro}
\epsilon_{\underline{d}}(i)=(-1)^{p_i a_{<i}+q_i (a_{<i}+1)+\frac{a_i(a_i-1)}{2}}.
\end{equation}
Le groupe $A(\psi_G)$ s'identifie de manière naturelle à   $(\pm 1)^\ell$ et 
$\epsilon_{\underline{d}}$ à un caractère de $A(\psi)$. 
Il découle alors essentiellement des résultats de \cite{AdJo87}, \cite{Joh1} et de \cite{AMR}  (voir aussi \cite{MR3}) que 
\begin{equation} \label{piApsiD} \pi^A(\psi)=
\sum_{\underline{d}\in \caD(\psi)} \scrA_{\underline{d} }(\psi)\boxtimes \epsilon_{\underline{d}}.\end{equation}
D'après la remarque faite ci-dessus sur le fait que les  $\scrA_{\underline{d} }(\psi)$
sont non isomorphes deux à deux, on en déduit la propriété de multiplicité un pour ces paquets.

Ensuite, on abandonne l'hypothèse de régularité  (\ref{regintro})  pour ne conserver que 
l'hypothèse de décroissance de la suite $(t_i)_{i=1,\ldots,\ell}$, et l'on définit $\scrA_{\underline d}(\psi)$ comme ci-dessus. L'induction cohomologique
a alors lieu dans le \og weakly fair range\fg\, et l'on a toujours 
 $\left(\caR_{\frqqq_{\underline d}, L_{\underline d}\cap K}^{\frg,K}\right)^{k} (\Lambda_{\underline d})=0$,
 si  $k \neq \dim(\frv_{\underline d}\cap \frk)$. De plus,   si $\scrA_{\underline d}(\psi)$ n'est pas nul, c'est un module unitaire
et irréductible, cette dernière propriété étant propre aux groupes unitaires  ({\sl cf.} \cite{Mat96}, \cite{Trap}) et est due à Barbasch et Vogan. 
Le groupe $A(\psi_G)$ s'identifie maintenant  à   un quotient de $\{\pm 1\}^\ell$. Nous établissons que le caractère $ \epsilon_{\underline{d}}$
défini  ci-dessus ce factorise par ce quotient si $\scrA_{\underline{d} }(\psi)$ est non nul (Proposition  \ref{propker}).
Nous montrons que le terme de droite de  (\ref{piApsiD}), qui est donc encore  bien défini comme représentation de $G\times A(\psi_G)$, est 
bien la représentation $\pi^A(\psi)$. On a donc (Théorème \ref{ThmApsiD} du texte)
\begin{thm} On suppose que $\psi$ est de bonne parité. Alors la représentation  associée à $\psi$ est 
\[\pi^A(\psi)=
\sum_{\underline{d}\in \caD(\psi)} \scrA_{\underline{d}}(\psi)\boxtimes \epsilon_{\underline{d}}.\]
De plus, les représentations $\scrA_{\underline{d}}(\psi)$ non nulles   sont non isomorphes deux à deux.
\end{thm}
Comme dans le cas des groupes classiques on passe du cas régulier au cas général (de bonne parité) en utilisant les foncteurs de translation.
Pour les groupes unitaires, on a donc en plus le fait que les $\scrA_{\underline d}(\psi)$ sont irréductibles ou nuls (on ne détermine pas quand ces modules sont nuls, 
voir la remarque \ref{rmqthm}).
De plus, la seconde assertion du théorème montre  que la propriété de multiplicité un est conservée.

Passons maintenant à un paramètre  $\psi_G$ général, avec $\psi=\psi_{mp}\oplus \psi_{bp}$.  A la partie de bonne parité $\psi_{bp}$ on attache donc un paquet $\Pi(\psi_{G_{bp}})$ 
d'un groupe unitaire de rang  $N_{bp}$ plus petit par la construction que l'on vient de donner, à  la partie de mauvaise parité
on attache une représentation irréductible unitaire $\rho$ d'un groupe $\GL_{N_\rho}(\bbC)$, et l'on a $N=2N_\rho+N_{bp}$.
Ceci détermine un sous-groupe parabolique standard $P=MN$ de $G$,
du moins si $\inf(p,q)\geq N_\rho$,  avec un facteur de Levi $M$ isomorphe à $GL_{N_\rho}(\bbC) \times 
G_{bp}$ (et donc $G_{bp}\simeq \U(p-N_\rho,q-N_\rho)$). On a alors (Théorème \ref{ThmReduc} du texte), en remarquant que les 
groupes $A(\psi_G)$ et $A(\psi_{G_{bp}})$ sont naturellement isomorphes,
\begin{thm} La représentation $\pi^A(\psi_G)$ est obtenue à partir de $\pi^A(\psi_{G_{bp}})$ par induction parabolique, c'est-à-dire
\[  \pi^A(\psi_G)=\bigoplus_{\pi_{bp}\in \Pi(\psi_{G_{bp}})}  \Ind_P^G(\rho \boxtimes \pi_{bp}) \boxtimes \rho_{\pi_{bp}}. \]
De plus, les induites paraboliques dans le membre de droite sont irréductibles.
\end{thm}

Dans cette étape, la  conservation de la propriété de multiplicité un des paquets est évidente. On a donc 
\begin{thm} Les paquets d'Arthur des groupes unitaires réels ont la propriété de multiplicité un.  
\end{thm}

\noindent {\sl Remerciements}. Nous remercions P. Trapa qui nous a indiqué que l'on pouvait déduire le résultat de multiplicité un du théorème 1.1 
d'une lecture attentive des résultats de la littérature sur les propriétés des foncteurs de translation.
Le deuxième auteur a bénéficié d'une aide de  l'agence nationale de la recherche 
ANR-13-BS01-0012 FERPLAY.

\section{Décomposition des $A$-paramètres}\label{decAparU}

Supposons d'abord que  $\mathbf G$ est un groupe classique. On note $\Std_{\mathrm G}$ la représentation standard du $L$-groupe de $\mathbf G$
dans $\GL_N(\bbC)$ (voir \cite{MR3}, \S 3.1), par exemple si $\mathbf G=\Sp_{2n}(\bbR)$, ${}^LG=\SO_{2n+1}(\bbC)\times W_\bbR$ et 
$ \Std_{\mathrm G}$ est donné par l'inclusion de $\SO_{2n+1}(\bbC)$ dans $\GL_{2n+1}(\bbC)$.

Si $\psi_G: \, W_\bbR \times \SL_2(\bbC)\rightarrow {}^LG$ est un paramètre d'Arthur pour $G$, on  pose
$\psi=\Std_G\circ \psi_G$ et l'on voit $\psi$ comme une représentation complètement réductible  de $W_\bbR\times \SL_2(\bbC)$.
Dans \cite{MR3}, \S 4.1, on a décomposé  cette  représentation en représentations irréductibles,
et séparé ces représentations irréductibles selon leur {\sl parité}, qui peut être {\sl bonne ou mauvaise},   ce qui permet 
 d'énoncer certains résultats de réduction. 

Nous allons maintenant faire de même pour les groupes unitaires.
Soit donc maintenant $\mathbf G$ un groupe unitaire réel de rang $N$, et soit 
\begin{equation} \label{AparU}\psi_G:  \, W_\bbR \times \SL_2(\bbC)\longrightarrow {}^LG  =\GL_N(\bbC)\rtimes W_\bbR\end{equation}  un paramètre  d'Arthur pour $G$.
On note $\psi$ la restriction de $\psi_G$ 
au sous-groupe $\bbC^\times \times \SL_2(\bbC)$ de $W_\bbR \times \SL_2(\bbC)$,
 que l'on voit comme une représentation de $  \bbC^\times \times \SL_2(\bbC)$ de dimension $N$. Cette représentation est complètement réductible. 
Pour tout $a\in \bbZ_{>0}$, notons
$R[a]$ la représentation algébrique de $\SL_2(\bbC)$ de dimension $a$, et pour tout 
$(t,s)\in \bbZ\times i\bbR$, notons
\begin{equation}  \label{chitx} \chi_{t,s}: \, z\mapsto (z/\overline{z})^{t/2} (z\overline{z})^{s/2}= z^{\frac{t+s}{2}} \, \bar z^{\frac{-t+s}{2}}  . \end{equation}
C'est un caractère unitaire  de $\bbC^\times$. On note pour tout $a\in \bbZ_{>0}$, 
\begin{equation}  \label{chitsa} \chi_{t,s,a}=\chi_{t,s}\circ {\textstyle \det_a} \end{equation}
où $ \det_a$ est le déterminant de $\GL_a(\bbC)$.

La forme générale de la  décomposition de $\psi$ en irréductibles  (après une éventuelle conjugaison dans $\GL_N(\bbC)$)   est 
\begin{equation}  \label{decApar} 
\psi=\bigoplus_{(t,s,a)\in \mathcal{E}(\psi)}\chi_{t,s}\otimes R[a] \end{equation}
pour un certain ensemble avec multiplicités finies $\mathcal{E}(\psi)$  de triplets  $(t,s,a)\in \bbZ\times i\bbR\times \bbZ_{>0}$.

\begin{defi}
On dit que le triplet $(t,s,a)$ est de bonne parité si $s=0$, et si $\frac{t+a-1}{2}+ \frac{N-1}{2}\in \mathbb{Z}$.
\end{defi}

\begin{rmq} Le fait que $\psi$ provient  d'un $A$-paramètre  pour $G=\U(p,q)$ est équivalent   a la propriété suivante  : 
si $(t,s,a)\in \mathcal{E}(\psi)$ n'est pas de bonne parité, alors la multiplicité de ce triplet dans $\mathcal{E}(\psi)$ est égale à  la multiplicité du triplet $(t,-s,a)$ 
dans le cas où  $s\neq 0$, et si $s=0$, la multiplicité de $(t,0,a)$ est paire. 
Cela revient à dire, en notant $\mathcal{E}(\psi)_{bp}$ les triplets de $\mathcal{E}(\psi)$ ayant bonne parité, 
qu'il existe  une décomposition de $\mathcal{E}(\psi)$ en l'union disjointe  de trois sous-ensembles, $\mathcal{E}(\psi)_{bp}$,
 $\mathcal{E}'(\psi)$,  et $\mathcal{E}''(\psi)$ tel que si $(t,s,a)\in \mathcal{E}'(\psi)$ alors $(t,-s,a)\in \mathcal{E}''(\psi)$ avec la même multiplicité.
\end{rmq}

On note 
\begin{equation}  \label{psibp} 
\psi_{bp}:= \bigoplus_{(t,s,a)\in \mathcal{E}(\psi)_{bp}} \chi_{t,s}\otimes R[a].\end{equation}
 Alors $\psi_{bp}$ est un morphisme de 
$\mathbb{C}^\times \times \SL_2(\bbC)$ dans $\GL_{N_{bp}}(\bbC)$ où 
\begin{equation}\label{nbp}
N_{bp}=\sum_{(t,s,a)\in\mathcal{E}(\psi)_{bp}}a,\end{equation}
 qui provient d'un $A$-paramètre comme en (\ref{AparU}), mais pour les  groupes unitaires de  rang $N_{bp}$.

\section{Réalisation des groupes unitaires et de leur $c$-Levi. Induction cohomologique}

\subsection{Paires paraboliques}\label{parthetastable}
Soit  $G$  le groupe des points réels d'un groupe algébrique connexe réductif défini sur $\bbR$. On fixe une involution de 
Cartan $\theta$ de $G$, et l'on note $K$ le sous-groupe des points fixes de $\theta$ : c'est un sous-groupe compact maximal de $G$.
On suppose que  $G$ et $K$ sont de même rang;    autrement dit, $G$ possède un sous-groupe de Cartan $T$
inclus dans $K$ et donc compact.
On note $\frt_0$, $\frk_0$  et $\frg_0$ les  algèbres de Lie respectives de $T$,$K$  et $G$ et $\frt$, $\frk$  et $\frg$ leur complexifiées.
Les  sous-algèbres paraboliques $\theta$-stables de $\frg$ sont obtenues de la manière suivante. On fixe un 
élément $\nu\in \sqrt{-1} \frt_0^*$, et l'on pose :
\begin{equation} \label{qnu}  \frl=\frg^\nu=\frt\oplus \left( \bigoplus_{\alpha\in \Delta(\frg,\frt), \bil{\nu}{\alpha}=0}\frg^\alpha\right), 
\; \frv =\bigoplus_{\alpha\in \Delta(\frg,\frt), \\ \bil{\nu}{\alpha}>0}\frg^\alpha, \; \frqqq=\frl\oplus \frv,  \;  L=\mathrm{Norm}_G(\frqqq).
\end{equation}
Dans cet article, nous appellerons paire parabolique 
une paire $(\frqqq,L)$ obtenue comme ci-dessus, avec $\frqqq$ sous-algèbre parabolique $\theta$-stable de $\frg$. Le  
 sous-groupe $L$ de $G$  sera appelé $c$-Levi de $G$ (terminologie de Shelstad \cite{She15}).

\medskip 

On note $\left(\caR_{\frqqq,L\cap K}^{\frg,K}\right)^k$ le foncteur d'induction cohomologique de Vogan-Zuckerman ({\sl cf.} \cite{Vgreen}, \S 6.3.1)
en degré $k$, de la catégorie des $(\frl,K\cap L)$-modules vers la catégorie des $(\frg,K)$-modules.
Dans ce contexte, le degré qui nous intéresse particulièrement, et même exclusivement,  est $S=\dim(\frv \cap \frk)$, et dans l'article, nous écrirons
$\left(\caR_{\frqqq,L\cap K}^{\frg,K}\right)^S$ sans préciser de nouveau ce qu'est $S$.

Si $\Lambda$ est un caractère unitaire  de $L$, on note  $\lambda$ sa différentielle, que l'on voit comme un élément
 de $i\frt_0^*$.
On pose alors 
\[A_\frqqq(\Lambda)=\left(\caR_{\frqqq,L\cap K}^{\frg,K}\right)^S (\Lambda).\]
Si le groupe $L$ est  connexe, $\lambda$ détermine $\Lambda$ et l'on note alors cette représentation $A_\frqqq(\lambda)$.

Sous certaines conditions sur le caractère infinitésimal de la représentation $\sigma$ de $L$ que l'on induit, on a des résultats d'annulation, d'irréductibilité et d'unitarité
des modules $\left(\caR_{\frqqq,L\cap K}^{\frg,K}\right)^S(\sigma)$. Nous renvoyons à \cite{KnVo} pour les définitions du (weakly) good range, du  (weakly) fair  range
et des représentations faiblement unipotentes, pour lesquels on a les résultats
suivants : dans le weakly good range  $\left(\caR_{\frqqq,L\cap K}^{\frg,K}\right)^k(\sigma)$ est nul si $k\neq S=\dim(\frv \cap \frk)$.
Si $\sigma$ est irréductible et dans le good range  (resp. weakly good range), $\caR_{\frqqq,L,G}^S(\sigma)$ est irréductible (resp. irréductible ou nul).
Si $\sigma$ est unitaire  et dans le weakly good range, $\left(\caR_{\frqqq,L\cap K}^{\frg,K}\right)^S(\sigma)$ est unitaire.
Si $\sigma$ est faiblement unipotente et dans le weakly fair range, alors $\left(\caR_{\frqqq,L\cap K}^{\frg,K}\right)^k(\sigma)$ est nul si $k\neq S$ et 
$\left(\caR_{\frqqq,L\cap K}^{\frg,K}\right)^S(\sigma)$ est unitaire si $\sigma$ est de plus unitaire. 
En revanche, on n'a pas de résultat d'irréductibilité en général dans le weakly fair range, ni même le fair range.

\subsection{Le groupe symplectique}\label{symp}
Une manière commode de réaliser les groupes unitaires $\U(p,q)$, $p+q=N$, est de les réaliser comme $c$-Levi
d'un groupe symplectique $\Sp(2N,\bbR)$, les $c$-Levi des $\U(p,q)$
s'identifient alors eux aussi à des $c$-Levi de $\Sp(2N,\bbR)$.
On suppose $\bbR^{2N}$ (identifié à $\caM_{2N,1}(\bbR)$, les matrices colonnes) muni de sa forme symplectique usuelle, c'est-à-dire, si $X,Y\in \bbR^{2N}$ 
\[ (X\vert Y)={}^tX J Y \text{ où }  J=\begin{pmatrix}0_N&I_N\\-I_N&0_N\end{pmatrix}.\]
Soit $\scrG=\Sp(2N,\bbR)$ le groupe des isomorphismes de $(\bbR^{2N},  (.\vert .))$, que l'on muni
de l'involution de Cartan $\theta: g\mapsto{}^tg^{-1}$. Le sous-groupe de $\scrG$ des  points fixes sous $\theta$
est un sous-groupe compact maximal de $\scrG$, que l'on note $\scrK$, et qui est isomorphe
au groupe unitaire $\U(N)$. On note $\frg_0$ et $\frk_0$ les sous-algèbres de Lie respectives de 
$\scrG$ et $\scrK$, réalisées comme sous-algèbre de Lie de $\caM_{2N}(\bbR)$.
Pour tout $(a_1,\ldots ,a_N)\in \bbR^N$, on pose :
\[t(a_1,\ldots,a_N)=\left(  \begin{array}{llll|llll}   & & & & a_1&&&\\
& & && & a_2 &&\\
&&& &&&\ddots&\\
&&& &&&& a_N\\
\hline
-a_1 & && & & &&\\
       &-a_2 &&& &&&\\
    & &\ddots &&   &&&\\
    &&& -a_N&&&&
  \end{array}\right). \]
Alors $\frt_0:=\{ t(a_1,\ldots,a_N), \, (a_1,\ldots ,a_N)\in \bbR^N \}$ est une sous-algèbre de Cartan de $\frk_0$ 
et aussi de $\frg_0$.

Notons $\frg,\frk,\frt$ les complexifications des algèbres de Lie $\frg_0,\frk_0,\frt_0$,
respectivement. Soient $\Delta(\frg,\frt)$, $\Delta(\frk,\frt)$ les systèmes de racines de 
$\frg$ et $\frk$ respectivement, relativement à la sous-algèbre de Cartan $\frt$.
On a 
\[\Delta(\frg,\frt)=\{ \pm(e_i \pm e_j), \, 1\leq i<j \leq N\} \cup \{ \pm 2e_i,  1\leq i \leq N\} , \]
\[\Delta(\frk,\frt)=\{ \pm(e_i - e_j), \, 1\leq i<j \leq N\}, \]
où $e_i\in  \sqrt{-1}\,  \frt_0^*\subset \frt^*$ est la forme linéaire $t(a_1,\ldots,a_N)\mapsto   \sqrt{-1} \, a_i$.
On fixe les systèmes de racines positives 
\[\Delta^+(\frg,\frt)=\{ (e_i \pm e_j), \, 1\leq i<j \leq N\} \cup \{  2e_i,  1\leq i \leq N\} , \]
\[\Delta^+(\frk,\frt)=\{ (e_i - e_j), \, 1\leq i<j \leq N\}. \]
On identifie $\frt^*$ et $\bbC^N$ grâce à la base $(e_i)_{1\leq i\leq  N}$ de $\frt^*$, et de même pour $\frt$ grâce à la base duale.

\subsection{Paires paraboliques maximales et induction cohomologique}\label{ppm}

On continue avec les notations du paragraphe précédent, en particulier $\scrG=\Sp(2N,\bbR)$.
Soit $(\frqqq,L)$ une paire parabolique pour $\scrG$, et l'on suppose que la sous-algèbre parabolique $\theta$-stable $\frqqq=\frl\oplus\frv$ est 
maximale.  Une telle sous-algèbre est obtenue en prenant un élément de $\frt$ de la forme
\[  t_{p,q}=(\underbrace{1,\ldots,1}_{p}, \underbrace{0,\ldots,0}_{N-p-q},  \underbrace{-1,\ldots,-1}_{q}) \]
avec $p+q\leq N$.
On pose alors 
\begin{equation}\label{qpq} \frl_{p,q}=\frg^{t_{p,q}}=\frt\oplus\left(  \bigoplus_{\alpha\in \Delta(\frg,\frt)\vert \alpha(t)=0 } \frg_\alpha  \right), 
\quad \frv'_{p,q}= \bigoplus_{\alpha\in \Delta(\frg,\frt)\vert \alpha(t)>0 } \frg_\alpha, \quad \frqqq'_{p,q}=\frl_{p,q}\oplus\frv'_{p,q}. 
\end{equation}
\begin{equation}\label{Lpq}L_{p,q}=\mathrm{Norm}_\scrG(\frqqq_{p,q})
\end{equation}

Dans le cas où $p+q=N$, $L_{p,q}$ est une réalisation du  groupe unitaire  $\U(p,q)$.
Les racines de $\frt$ dans $\frl_{p,q}$ sont :
\begin{align*}
&\pm(e_i-e_j), &\,  1\leq i<j\leq p  \text{ ou }   p+1 \leq i<j\leq N,\\
&\pm(e_i+e_j),& \,   1\leq i \leq p< j\leq N. 
 \end{align*}
 On choisit comme système de racines positives :
 \[\Delta_{p,q}^+=\Delta^+(\frl_{p,q},\frt)= \left\{  \begin{array}{ll}  e_i-e_j,  &  1\leq i<j\leq p, \\
   -(e_i-e_j) , &  p+1 \leq i<j\leq N \\
  e_i+e_j, & \,   1\leq i \leq p< j\leq N
    \end{array} \right\}.\]
 On pose 
 \begin{equation}\label{deltalpq}
 \delta(\frl_{p,q})=\frac{1}{2}\! \! \! \sum_{\alpha\in \Delta^+_{p,q}}   \! \! \! \alpha
  ={\small \frac{1}{2} (\underbrace{N-1,N-3 , \ldots, q-p+1}_{p} ,\underbrace{ p-q+1,\ldots,  p-q+3,\ldots ,N-1}_{q}  )}. \end{equation}

On a $\scrK=L_{N,0}$ et $L_{p,q}\cap \scrK\simeq \U(p)\times \U(q)$.

\subsection{$c$-Levi des $\U(p,q)$}\label{cLUpq}
On continue avec les notations de la section précédente. Notons $\caD(p,q)$ l'ensemble des familles
$\underline d=(p_i,q_i)_{i=1,\ldots,\ell}$ de couples d'entiers positifs ou nuls tels que  $\sum_i p_i=p$ et $\sum_i q_i=q$ et 
$\caD(N)=\coprod_{p,q\vert p+q=N} \caD(p,q)$.
Pour tout $\underline d=(p_i,q_i)_{i=1,\ldots,\ell} \in \caD(N)$, notons 
\[  t_{\underline d}=(\underbrace{\ell,\ldots,\ell}_{p_1},\ldots, \underbrace{2,\ldots,2}_{p_{\ell-1}},   \underbrace{1,\ldots,1}_{p_\ell}, 
 \underbrace{-1,\ldots,-1}_{q_\ell}, \underbrace{-2,\ldots,-2}_{q_{\ell-1}}, \ldots , \underbrace{-\ell,\ldots,-\ell}_{q_1}) .\]
On définit   la sous-algèbre parabolique $\theta$-stable $\frqqq'_{\underline d}=\frl_{\underline d}\oplus \frv'_{\underline d}$ et le $c$-Levi 
$L_{\underline d}$ associés à $t_{\underline d}$ comme en (\ref{qnu}).
Si $\underline d\in \caD(p,q)$, alors $\frl_{\underline d}\subset \frl_{p,q}$ et  $L_{\underline d}$ est un $c$-Levi de $L_{p,q}$,
associé à la sous-algèbre parabolique $\theta$-stable $\frqqq_{\underline d}$ de $\frl_{p,q}$, où
\begin{equation}\label{qd}  
\frqqq_{\underline d}= \frqqq'_{\underline d}\cap \frl_{p,q}=   \frl_{\underline d}  \oplus  (\frl_{p,q}\cap \frv'_{\underline d})
=\frl_{\underline d}  \oplus   \frv_{\underline d}. 
\end{equation}
On pose $\Delta^+(\frl_{\underline d},\frt)=\Delta(\frl_{\underline d},\frt)\cap \Delta^+(\frl_{p,q},\frt)$ et pour tout $i=1,\ldots,\ell$,$a_i=p_i+q_i$
 \begin{equation} \label{deltald}\delta(\frl_{\underline d})
 =\frac{1}{2}\sum_{\alpha\in \Delta^+(\frl_{\underline d},\frt)}\alpha\end{equation}
{\tiny \[  =\frac{1}{2} (\underbrace{a_1-1,a_1-3 , \ldots, q_1-p_1+1}_{p_1} ,\underbrace{a_2-1,a_2-3 , \ldots, q_2-p_2+1}_{p_2}, \ldots , 
\underbrace{a_\ell-1,a_\ell-3 , \ldots, q_\ell-p_\ell+1}_{p_\ell}, \]
  \[\underbrace{ p_\ell-q_\ell+1,\ldots,  p_\ell-q_\ell+3,\ldots ,a_\ell-1}_{q_\ell},\ldots,
 \underbrace{ p_2-q_2+1,\ldots,  p_2-q_2+3,\ldots ,a_2-1}_{q_2},\underbrace{ p_1-q_1+1,\ldots,  p_1-q_1+3,\ldots ,a_1-1}_{q_1}  ). \]
}

\section{Paramètres de bonne parité} \label{bonneparite} 
Soit $\psi$ un $A$-paramètre pour $\U(p,q)$ comme en  (\ref{AparU}).
On suppose que $\psi$ est de bonne parité, {\sl i.e.}  $\psi=\psi_{bp}$. On a donc 
\begin{equation}  \label{decAparbp} 
\psi=\bigoplus_{i=1}^\ell (\chi_{t_i}\otimes R[a_i]) .\end{equation}
On a supprimé les paramètres $s_i$ des notations, puisque dans le cas de bonne parité, ils sont tous nuls.
 On ordonne les indices $i$ pour que 
   $t_1\geq \cdots \geq t_\ell$ et  si $t_i=t_{i+1}$ alors $a_i\geq a_{i+1}$.

On reprend les notations de la section \ref{cLUpq}, en particulier la réalisation de  $\U(p,q)$ comme  sous-groupe
$L_{p,q}$ d'un groupe symplectique.
Pour tout $i\in [1,\ell]$, on fixe une décomposition $a_i=p_i+q_i$ en somme de deux entiers éventuellement nuls. On note $\mathcal{D}(\psi)$ 
l'ensemble de ces décompositions qui vérifient  de plus  $p=\sum_i p_i$. 
On a donc  $\mathcal{D}(\psi)\subset \caD(p,q)$, ce dernier ensemble étant défini dans la section    \ref{cLUpq}, et  en particulier 
un élément   $\underline{d}\in \mathcal{D}(\psi)$
détermine un $c$-Levi  $L_{\underline d}$  de $L_{p,q}$ isomorphe à $ \times_i \U(p_i,q_i)$ et 
  une sous-algèbre   parabolique  $\theta$-stable $\frqqq_{\underline d}=\frl_{\underline d}\oplus \frv_{\underline d}$ de 
  $\frl_{p,q}$ dont la sous-algèbre de Levi est  précisément l'algèbre de Lie complexifiée de $L_{\underline{d}}$. 
Pour tout $i\in [1,\ell]$, on pose $a_{<i}:=\sum_{j<i}a_j$.

On définit  le  caractère $\Lambda_{\underline{d}}$ du groupe $L_{\underline{d}}$ via l'isomorphisme de $L_{\underline{d}}$ avec $ \times_i \U(p_i,q_i)$
par 
\begin{equation}\label{lambdads}
\Lambda_{\underline d}=\boxtimes_i  \textstyle \det^{\frac{t_i+a_i-N}{2}-a_{<i} }.
\end{equation}
On peut aussi définir ce caractère en donnant sa différentielle comme un élément de $\frt^*$ dans le système de coordonnées
de la section \ref{symp}. Posons $\lambda_i=\frac{t_i+a_i-N}{2}-a_{<i} $. On a alors
\[ d\Lambda_{\underline d}={
\left(\underbrace{\lambda_1,\ldots ,\lambda_1}_{p_1}, \ldots , \underbrace{\lambda_\ell,\ldots ,\lambda_\ell }_{p_\ell},
\underbrace{-\lambda_\ell,-\ldots ,\lambda_\ell }_{q_\ell},\ldots, \underbrace{-\lambda_1,\ldots,- \lambda_1 }_{q_1} \right)} \]
 Le foncteur d'induction cohomologique   $\left(\caR_{\frqqq_{\underline d}, L_{\underline d}\cap K}^{\frg,K}\right)^{\dim(\frv_{\underline d}\cap \frk)}$
   permet de produire une représentation de $\U(p,q)$ à partir de ce caractère.
On pose  
\begin{equation}\label{lqd}
\scrA_{\underline d}(\psi)=\left(\caR_{\frqqq_{\underline d}, L_{\underline d}\cap K}^{\frg,K}\right)^{\dim(\frv_{\underline d}\cap \frk)} (\Lambda_{\underline d}).
\end{equation}
La condition  $t_1\geq \cdots \geq t_{\ell}$  assure que cette induction cohomologique est dans le weakly fair range ({\sl cf.} \cite{KnVo} p. 35). 
De ceci, il découle que $\left(\caR_{\frqqq_{\underline d}, L_{\underline d}\cap K}^{\frg,K}\right)^{i} (\Lambda_{\underline d})=0$,
 si  $i\neq \dim(\frv_{\underline d}\cap \frk)$, et que  si $\scrA_{\underline d}(\psi)$ n'est pas nul, c'est un module unitaire
et irréductible, cette dernière propriété étant propre aux groupes unitaires  ({\sl cf.} \cite{Mat96}, \cite{Trap}).

On  a défini  l'application $\epsilon_{\underline{d}}$ de $[1,\ell]$ dans $\pm 1$ en (\ref{edintro}).
Nous allons expliquer plus loin comment $A(\psi)$ s'identifie à un quotient de  $\{\pm 1\}^\ell$ et 
$\epsilon_{\underline{d}}$ à un caractère de $A(\psi)$. 

\begin{thm}  \label{ThmApsiD}On suppose que $\psi$ est de bonne parité. Alors la représentation associée à $\psi$ est 
\[\pi^A(\psi):=
\sum_{\underline{d}} \scrA_{\underline{d}}(\psi)\boxtimes \epsilon_{\underline{d}}.\]
De plus, les représentations $\scrA_{\underline{d}}(\psi)$ non nulles sont non isomorphes deux à deux.
\end{thm}

\ 
\begin{rmq} \label{rmqthm} On ne précise pas ici  quand les $\scrA_{\underline{d}}(\psi)$ sont non nulles.
 Toutefois l'article \cite{Trap} donne un algorithme pour résoudre cette question;  c'est un problème difficile.
 La deuxième assertion du théorème est une assertion de multiplicité un dans les paquets d'Arthur pour des paramètres de bonne parité.
 Les résultats de réduction à la bonne parité établis dans la section \ref{irredU}
  montrent qu'il y a multiplicité un dans tous les paquets d'Arthur des groupes unitaires.
\end{rmq}
\ 

\dem  Le théorème est  déjà établi  sous la condition que les induites  cohomologiques (\ref{lqd}) soient réalisées dans le good range, 
ce qui est plus restrictif que le weakly fair range, la condition sur $\psi$ étant (\ref{regintro}).
 Dans ce cas, $A(\psi_G)$ s'identifie naturellement à $\{\pm 1\}^\ell$ . 
  
En effet, sous cette condition, le théorème est un cas particulier de  \cite{MR3}, Théorème 9.3,  qui 
 résulte  fondamentalement de  \cite{AdJo87},  \cite{Joh1} et  \cite{AMR}.
De plus, dans ce cas, les  $\scrA_{\underline{d}}(\psi)$ sont non nulles et non isomorphes deux à deux.

On  se ramène à ce cas en utilisant les foncteurs de translation ({\sl cf.} \cite{KnVo} chapter VII et chapter VII, \S5).
On fixe  des entiers $T_1>> \cdots>> T_\ell\geq 0$ et on note $\psi_+$ le $A$-paramètre
\begin{equation}  \label{psi+} 
\psi_+=\bigoplus_{i} (\chi_{t_i+T_i}\otimes R[a_i]),  \end{equation}
de sorte que $\psi_+$ vérifie les hypothèse de good range (\ref{regintro}), et l'on a donc 
\[
\pi^A(\psi_+)=\bigoplus_{\underline{d}\in \caD(\psi_+)} \scrA_{\underline{d}}(\psi_+)\boxtimes \epsilon_{\underline{d}},
\]
où les modules  $\scrA_{\underline{d}}(\psi_+)$ sont définis comme en (\ref{lqd}) en remplaçant les 
$t_i$ par les $t_i+T_i$ dans la définition de $\Lambda_{\underline d}$ en (\ref{lambdads}).
Avec ces hypothèses, comme nous l'avons dit, $A(\psi_+)$ s'identifie  à $\{\pm 1\}^\ell$  et cette expression est alors bien définie.

On considère  le foncteur de translation correspondant à la représentation de dimension finie $\scrF$ qui est la restriction à $\U(p,q)$
 de la représentation de $\GL_N(\mathbb{C})$ de plus bas poids (dans les coordonnées usuelles pour ce groupe)
  \[ \left( \underbrace{-T_1/2, \cdots, -T_1/2}_{a_1}, \cdots, \underbrace{-T_\ell/2, \cdots, -T_\ell/2}_{a_\ell}\right). \]
 Dans les coordonnées choisies  ci-dessus pour $\frt^*$, ce plus bas poids est 
 \[ \gamma=
\left(\underbrace{-T_1/2,\ldots ,-T_1/2}_{p_1}, \ldots , \underbrace{-T_\ell/2,\ldots ,-T_\ell/2 }_{p_\ell},
\underbrace{T_\ell/2,-\ldots ,T_\ell/2 }_{q_\ell}, \ldots,  \underbrace{T_1/2,\ldots,T_1/2 }_{q_1} \right). \]

\begin{rmq}\label{mdeltau}
 Posons $\delta(\frv_{\underline d})=\delta(\frl_{p,q})-\delta(\frl_{\underline d})$. On voit facilement avec les expressions (\ref{deltalpq}) 
et (\ref{deltald}) que l'on peut définir les $T_i$ vérifiant les conditions ci-dessus en posant  $\gamma=-2m\, \delta(\frv_{\underline d})$
 pour un entier $m$ assez grand.
 \end{rmq}
 
  On note $\mathcal{T}$ ce foncteur et on sait d'après \cite{MR5}, théorème 4.1  que l'on a 
  \[
\mathcal{T}(\pi^A(\psi_+))=\pi^A(\psi),
\]
et d'autre part, d'après \cite{MR5}, théorème 4.3, en  remarquant que  $\caD(\psi_+)=\caD(\psi)$, 
on a $ \mathcal{T}( \scrA_{\underline{d}}(\psi_+))= \scrA_{\underline{d}}(\psi)$ pour tout 
$\underline d\in \caD(\psi)$. Pour un tel 
 $\underline{d}$ fixé, le caractère $\epsilon_{\underline{d}}$ se factorise en un caractère de $A(\psi)$ si et seulement si pour tout $i,j\in [1,\ell]$ 
 tels que $t_i=t_j$ et $a_i=a_j$, on a $\epsilon_{\underline{d}}(i)=\epsilon_{\underline{d}}(j)$.
Or on a d'après \cite{Mat96} et \cite{Trap} que l'on commentera ci-dessous

\begin{prop}\label{propker} La représentation $\mathcal{T}(\scrA_{\underline{d}}(\psi_+))$ est non nulle  seulement si pour tout $i\in [1,\ell[$ tel que 
$t_{i+1}=t_i$ on a $p_i\geq q_{i+1}$ et $q_i\geq p_{i+1}$. 
\end{prop}

Comme l'a remarqué Matumoto, cette proposition est une conséquence du paragraphe 4.2 de \cite{BVPIM}. 
On peut donc dire que cette proposition est essentiellement due à Barbasch et Vogan.
 On peut faire des inductions par étage et donc supposer que $\ell=2$,  $i=1$ et $t_1=t_2$. 
 On rappelle que l'on a ordonné les couples $(t_i,a_i)$ tel que si $t_i=t_{i+1}$ alors $a_i\geq a_{i+1}$ 
 sans cette hypothèse la proposition serait fausse. Cette hypothèse est un oubli dans \cite{Mat96} et les choix de cette référence 
 pour l'induction cohomologique étant opposés aux nôtres les inégalités sont inversées.

Ainsi pour tout $\underline d\in \caD(\psi)$, soit le caractère   $\epsilon_{\underline d}$ se factorise en un caractère de $A(\psi)$,  soit  
$\scrA_{\underline d}(\psi)$ est nul,  et ceci donne un sens
à la formule donnant $\pi^A(\psi)$ dans le théorème et établit la première partie de celui-ci.

 \ 

Démontrons maintenant l'assertion de multiplicité un. La difficulté vient de l'utilisation du foncteur de translation
permettant de passer du good range au weakly fair range.
Avec les notations ci-dessus, on sait que si $\underline d\neq \underline d'$, alors $\scrA_{\underline{d}}(\psi_+)$ et 
$\scrA_{\underline{d'}}(\psi_+)$ sont non isomorphes, et l'on voudrait en déduire que si les représentations 
$\scrA_{\underline{d}}(\psi) =\mathcal{T}(\scrA_{\underline{d}}(\psi_+))$ et $\scrA_{\underline{d'}}(\psi)=\mathcal{T}(\scrA_{\underline{d'}}(\psi_+))$
sont non nulles, alors elles sont non isomorphes. 

Pour cela, nous allons utiliser des résultats sur les foncteurs de translation disséminés dans \cite{KnVo}, et nous allons expliquer comment 
nous y ramener, ce qui nous oblige à quelques détours. Signalons aussi que les idées sont expliquées et mises en oeuvre dans \cite{VogDs}
dans un contexte différent, et que c'est P. Trapa qui nous a suggéré que la démonstration de l'irréductibilité de l'induction cohomologique
convenablement comprise pouvait aussi donner l'énoncé de multiplicité un voulu. 
La première chose à faire est de se ramener à des foncteurs d'induction cohomologique où la sous-algèbre parabolique
$\theta$-stable $\frqqq$ est fixée et où ce sont les formes réelles fortes au sens de \cite{ABV} qui vont varier.
Expliquons le formalisme. On part du groupe compact $\U(N)$ et de sa réalisation usuelle comme sous-groupe des points
fixes de l'involution $\sigma : g\mapsto {}^t \bar g^{-1}$ de $\GL_N(\bbC)$. On note ici $\frg$ l'algèbre de Lie de $\GL_N(\bbC)$.
On choisit un tore maximal $T$ de $\U(N)$, et l'on note $\frt$ la complexifiée de son algèbre de Lie.
On fixe une sous-algèbre de Borel $\frb$ de $\frg$ contenant $\frt$. On identifie  $T$ à $\U(1)^N$, $\frt$ à $\bbC^N$ de sorte que 
les racines simples de $\frt$ dans $\frb$ soient les formes linéaires $e_i-e_{i+1}$, $i=1,\ldots ,N-1$, où $(e_i)_{i=1,\ldots,N}$ est la base canonique de 
$(\bbC^N)^*$.  La partition $N=\sum_{i=1}^\ell a_i$ détermine alors
une sous-algèbre parabolique $\frqqq=\frl\oplus \frv$ de $\frg$ contenant $\frb$.
Notons $T[2]$ l'ensemble des éléments d'ordre 2 de $T$. Pour tout $t\in T[2]$, posons 
\[\sigma_t=\Ad(t)\circ \sigma, \quad \theta_t =\Ad(t).\]
Alors $\sigma_t$ est une forme réelle de $\GL_N(\bbC)$ et nous notons $\U_t$ le groupe de ses points réels.
L'involution $\theta_t$ est une involution de Cartan de $\U_t$, et l'on note $K_t$ le sous-groupe de ses points fixes. C'est un sous-groupe
compact maximal de $\U_t$. Les $\U_t$ sont des formes intérieures pures de $\U(N)=\U_{t=1}$. 
Si l'on écrit $t=(\eta_1,\ldots ,\eta_i)$ par l'identification  $T=\U(1)^N$, avec $\eta_i\in \{\pm 1\}$ pour tout $i=1,\ldots,N$,
alors $\U_t$ est isomorphe au groupe $\U(n_1(t),n_{-1}(t))$ où $n_1(t)$ est le nombre des $\eta_i$ égaux à $1$ et 
$n_{-1}(t)$ celui des $\eta_i$ égaux à $-1$. Pour $t_1,t_2\in T[2]$, les formes réelles $\U_{t_1}$ et $\U_{t_2}$ sont équivalentes
si et seulement si $n_{1}(t_1)=n_{1}(t_2)$ et $n_{-1}(t_1)=n_{-1}(t_2)$. 
La sous-algèbre parabolique $\frqqq$ est $\theta_t$-stable pour tout $t\in T[2]$. Posons 
$L_t=\mathrm{Norm}_{\U_t}(\frqqq)$ : c'est un $c$-Levi de $\U_t$ et la complexifiée de son algèbre de Lie est $\frl$.
Ces groupes $L_t$, pour $t\in T[2]$, sont des formes réelles d'un groupe complexe $L_\bbC$ isomorphe à $\times _i \GL_{a_i}(\bbC)$.
Le groupe $L_t$ est isomorphe à $\times_i \U(p_i,q_i)$, où 
$p_1$ (resp $q_1$) est le nombre de $1$ (resp. de $-1$) dans les $a_1$ premières coordonnées de $t$, et ainsi de suite.
On a donc pour chaque $t$ des foncteurs d'induction cohomologique
\[\left(\caR_{\frqqq,L_t\cap K_t}^{\frg,K_t}\right)^k \]
que nous allons toujours considérer dans le bon degré, c'est-à-dire $S_t=\dim (\frv \cap \frk_t)$.

On revient maintenant à $p$ et $q$ fixés avec $p+q=N$ et pour tout  élément $\underline d =(p_i,q_i) \in \caD(\psi)$, 
on pose 
\[t_{\underline d}=(\underbrace{1,\ldots,1}_{p_1},\underbrace{-1,\ldots,-1}_{q_1},\quad  \ldots \quad , \underbrace{1,\ldots,1}_{p_\ell},\underbrace{-1,\ldots,-1}_{q_\ell}).  \]

On fixe  un caractère $\Lambda$ de $L_\bbC$ : dans l'identification de $L_\bbC$ avec  $\times _i \GL_{a_i}(\bbC)$, il est donné par 
\[\Lambda=\boxtimes_i  \textstyle \det^{\frac{t_i+a_i-N}{2}-a_{<i} }. \]
On pose alors pour tout  $\underline d =(p_i,q_i) \in \caD(\psi)$,
\[ \scrA'_{\underline d}(\psi)= \left(\caR_{\frqqq,L_{t_{\underline d}} \cap K_{t_{\underline d}}}^{\frg,K_{t_{\underline d}}}\right)^{S_{t_{\underline d}}} (\Lambda).\]
C'est un $(\frg, K_{t_{\underline d}})$-module. Comme les formes réelles $\U_{t_{\underline d}}$ sont toutes équivalentes lorsque 
$\underline d$ décrit $\caD(\psi)$, c'est-à-dire conjuguées deux à deux par un automorphisme intérieur de $\GL_N(\bbC)$,
 on peut voir les classes d'équivalence de $(\frg, K_{t_{\underline d}})$-modules
comme des classes d'équivalence de modules de Harish-Chandra pour  $\U(p,q)$, et les modules  $\scrA'_{\underline d}(\psi)$
correspondent aux modules $ \scrA_{\underline d}(\psi)$. Il s'agit donc de montrer que les 
modules $\scrA'_{\underline d}(\psi)$ non nuls sont inéquivalents deux à deux, vus comme modules de Harish-Chandra pour $\U(p,q)$.
D'autre part,  comme les groupes unitaires
sont   connexes, la catégorie des  $(\frg, K_{t_{\underline d}})$-modules  est une sous-catégorie pleine de la catégorie des $\frU(\frg)$-modules, 
où $\frU(\frg)$ est l'algèbre enveloppante de $\frg$, et il s'agit donc de montrer que les 
modules $\scrA'_{\underline d}(\psi)$ non nuls sont inéquivalents deux à deux comme $\frU(\frg)$-modules.

Bien sûr, on peut remplacer $\psi$ par $\psi_+$ dans ces considérations, et l'on obtient de même des modules  $\scrA'_{\underline d}(\psi_+)$,   
$\underline d \in \caD(\psi_+)=\caD(\psi)$ et l'on sait qu'ils sont inéquivalents deux à deux comme $\frU(\frg)$-modules. 

Rappelons maintenant quelques éléments sur le foncteur de translation tirés de \cite{KnVo}, Chapter VII. Celui-ci est donc défini par la représentation
de dimension finie $\scrF$ de plus bas poids $\gamma=-2m \, \delta(\frv)$, d'après la remarque \ref{mdeltau}, et nous sommes donc dans les hypothèses
de la  Proposition 8.31 de \cite{KnVo}.
Le foncteur $\caT$ nous fait passer des modules ayant comme caractère infinitésimal généralisé celui donné par le paramètre $\psi_+$, notons-le $ \mu_+$,
à des modules ayant  comme caractère infinitésimal généralisé celui donné par le paramètre $\psi$, notons-le $ \mu=\mu_+ + \gamma=\mu_+ -2m\, \delta(\frv)$.
Ici, on voit $\mu_+$ et  $\mu$ comme des éléments de $\frt^*$ qui déterminent chacun  un caractère du centre de l'algèbre enveloppante.

Reprenons des éléments de la démonstration de la Proposition 8.31 de \cite{KnVo}, que l'on particularise au cas des groupes unitaires étudiés ici, en adaptant
légèrement les notations. On y introduit pour tout $\lambda'\in \frt^*$, un module de Verma généralisé
noté $M(\lambda')$, et l'on montre que si $\lambda'$ vérifie les conditions du weakly fair range, alors $\caT(M(\lambda'+2m\, \delta(\frv))=M(\lambda')$ (Lemma 8.35). 
On en déduit (Lemma 8.39) que pour toute forme réelle $G$ du groupe complexe ambiant muni d'une involution  de Cartan $\theta$ ayant  comme groupe des 
points fixes le sous-groupe compact maximal $K$ de $G$, tel que $\frqqq$ soit une sous-algèbre parabolique $\theta$-stable, 
 avec de plus  les conditions de weakly fair range  sur $\lambda'$,  que l'on a $\caT(A_{\frqqq,L\cap K}^{\frg,K}(\lambda'+2m\, \delta(\frv))=A_{\frqqq,L\cap K}^{\frg,K}(\lambda')$.
En particulier, ceci s'applique aux formes réelles $\U_{t_{\underline d}}$ définies ci-dessus et a $\lambda'+\delta=\mu$ 
($\delta$ désigne bien entendu la demi-somme des racines positive), et l'on a donc
\[ \caT(\scrA'_{\underline d}(\psi_+))= \scrA'_{\underline d}(\psi),\]
pour tout $\underline d\in \caD(\psi)$.
On introduit
 $Q(\lambda')=\End(M(\lambda'))$, qui est un $\frU(\frg)\otimes \frU(\frg)$-module ayant pour caractère infinitésimal $(\lambda'+\delta,-(\lambda'+\delta))$
muni d'une application naturelle $\varphi: \, \frU(\frg)\rightarrow Q(\lambda')$ qui respecte les actions à gauche de $\frU(\frg)$.
Ensuite, on prend les élément $\U(N)$-finis de l'algèbre $Q(\lambda')$ en posant $R(\lambda')=\End(M(\lambda'))_{\U(N)}$.
Ainsi, $R(\lambda')$ devient un $(\frg\oplus \frg, \U(N)\times \U(N))$-module, et l'image de $\phi$ est à valeurs dans $R(\lambda')$.
On  pose alors $S=S(\lambda')=R(\lambda')\otimes \End (\scrF)$, et cette algèbre admet une décomposition selon 
 ses composantes primaires (à gauche et à droite), $S=\oplus_{\alpha,\beta} S_\alpha^\beta$.

Prenons maintenant $\lambda'=\mu_+-\delta$.  La composante $S_\mu^{-\mu}$ est une sous-algèbre de $S$. 
Soit $M$ un $\frU(\frg)$-module ayant pour caractère infinitésimal $\mu_+$, et supposons que $M$ soit aussi un module à gauche unifère pour $R(\lambda')$
tel que les deux actions soient compatibles via  $\varphi$. Alors la composante $\mu$-primaire $N_\mu$  du $S$-module $N=M\otimes \scrF$ est naturellement 
un  $S_\mu^{-\mu}$-module. Ce qui est fondamental pour nous ici est le résultat suivant : si $M$ est un $R(\lambda')$-module simple, 
alors  $N_\mu$ est $S_\mu^{-\mu}$-module simple ou bien $0$, et de plus, si $M^1$, $M^2$ sont deux $R(\lambda')$-module simples non équivalents, et si 
 $N^1_\mu$ et $N^2_\mu$ sont non nuls, alors ce sont deux  $S_\mu^{-\mu}$-modules simples non équivalents.
On trouve la démonstration  à la page 524 de \cite{KnVo}, pour une algèbre $S$ différente, mais la démonstration est formelle et se transpose sans changements.
D'ailleurs un  résultat formel général analogue (avec la même démonstration) dans le cadre des algèbres à idempotents se trouve dans \cite{Ren10}, proposition I.3.2.

Il y a un foncteur de translation pour les $\frU(\frg)\otimes \frU(\frg)$-modules  construit avec $\End(\scrF)$ allant des modules ayant pour caractère infinitésimal 
généralisé $(\mu_+,-\mu_+)$ vers les modules ayant pour caractère infinitésimal généralisé $(\mu,-\mu)$. Notons le $\caT^2$. On a alors 
\[S_\mu^{-\mu}=\caT^2 (R(\mu_+-\delta))=R(\mu-\delta).\]
De plus, on peut munir les $\scrA'_{\underline d}(\psi_+)$, $\underline d\in \caD(\psi)$, d'une structure de $R(\lambda')$-module compatibles
avec l'action de $\frU(\frg)$, ceci apparaît à la page 577 de \cite{KnVo}.
Il résulte de ceci que les  $\scrA'_{\underline d}(\psi_+)$, $\underline d\in \caD(\psi)$ sont des $R(\mu-\delta)$-modules nuls ou simples, 
les non nuls étant inéquivalents deux à deux. La proposition 8.31 de \cite{KnVo}
a en fait pour but d'énoncer un critère pour en déduire que ce sont des $\frU(\frg)$-modules nuls ou simple, il suffit que l'application naturelle $\varphi$
de $\frU(\frg)$ dans $R(\mu-\delta)$ soit surjective. Or c'est le cas pour les groupes unitaires, cela vient du fait que les orbites nilpotentes en type A sont 
de Richardson avec une application moment birationnelle et  sont d'adhérence normale. C'est ainsi que l'on montre que les $A_\frqqq(\lambda)$-modules  
dans le weakly fair range 
sont nuls ou irréductibles pour les groupes unitaires. Ce que nous venons de remarquer ici, c'est que l'inspection
de la démonstration montre en plus que les modules $\scrA'_{\underline d}(\psi)$, $\underline d\in \caD(\psi)$ non nuls ne sont pas équivalents
en tant que $\frU(\frg)$-modules et 
ceci termine la démonstration du théorème.\qed

\begin{rmq}
La définition des paquets d'Arthur par  les identités de transfert endoscopiques suppose avoir choisi  parmi les formes
réelles fortes   (au sens de \cite{ABV}) $\U_t$, $t\in T[2]$,  introduites  ci-dessus, une forme quasi-déployée, et pour celle-ci, une donnée de Whittaker.
Ici, la forme quasi-déployée choisie est donné par 
\[t_*=(1,-1,1,-1,\ldots ,(-1)^{N-1})\]
 et est donc isomorphe à $\U(N/2,N/2)$ si $N$  est pair  et $\U(\lfloor N/2\rfloor+1,\lfloor N/2\rfloor)$ si $N$ est impair.
 D'autre part, on peut considérer un paramètre de séries discretes $\psi_G$,  pour les groupes unitaires de rang $N$,  c'est-à-dire  
avec $\psi=\oplus_{i=1}^N (\chi_{t_i}\boxtimes R[1])$ avec les $t_i$ distincts. Les constructions faites ci-dessus
pour $\underline d\in \caD(\psi)=((1,0),(0,1),(1,0),\ldots)$ déterminent une série discrète générique de ce groupe unitaire quasi-déployé, 
et l'on fixe la  donnée de Whittaker pour que cette série discrète admette une fonctionnelle de Whittaker.
Des choix différents mèneraient à une formule différente en (\ref{edintro}) en tordant la paramétrisation par un caractère de $A(\psi_G)$. 
\end{rmq}

\section{Réduction au cas de bonne parité}

\subsection{Mauvaise parité et induction parabolique}

Dans cette section, nous démontrons des résultats énoncés sans démonstration  dans \cite{MR3} ainsi que leurs analogues pour les groupes unitaires.
Soient $\mathbf G$ un groupe classique ou unitaire, et  $\psi_G$, $\psi$ comme dans la section \ref{decAparU}.
Considérons   une décomposition de $\psi$ de la forme :
\begin{equation}\label{psipsi}
       \psi=\rho\oplus \rho^*\oplus \psi' 
 \end{equation}
 où, dans $\rho$, il n'apparaît que des facteurs de mauvaise parité. Ici  $\rho^*$ désigne la représentation contragrédiente si 
 $\mathbf G$ est un groupe classique, et la duale hermitienne si $\mathbf G$ est un groupe unitaire.
 Remarquons que toute la partie de mauvaise parité peut se mettre sous la forme $\rho \oplus \rho^*$.

 Si $\mathbf G$ est classique,  le paramètre $\psi'$ se factorise par le $L$-groupe d'un groupe classique quasi-déployé 
$\mathbf {G}^\flat$ de même type que $\mathbf G$. Soit $\psi_{G^\flat}$ le paramètre d'Arthur pour le groupe $\mathbf G^\flat$ tel que
 $\psi'=\Std_{G'} \circ\psi_{G^\flat}$. De même, si $\mathbf G$ est unitaire, $\psi'$ est la restriction à $\bbC^\times  \times \SL_2(\bbC)$ d'un paramètre 
$\psi_{G'}$ pour un groupe unitaire quasi-déployé $\mathbf G^\flat$ de rang plus petit.

Notons $N_\rho$ la dimension de la représentation $\rho$.  Si $\mathbf G$ est classique, c'est une représentation 
de   $W_\bbR \times \SL_2(\bbC)$, et l'on note $ \Pi^{\GL}_\rho$
la représentation de $\GL_{N_\rho}(\bbR)$ de paramètre d'Arthur $\rho$ ({\sl cf. } \cite{AMR},   \S 3.1).
Si $\mathbf G$ est unitaire,  c'est une représentation  de $\bbC^\times \times   \SL_2(\bbC)$, 
et l'on note $ \Pi^{\GL}_\rho$
la représentation de $\GL_{N_\rho}(\bbC)$ de paramètre d'Arthur $\rho$. Pour unifier les notations, on 
note simplement $\GL_{N_\rho}$ pour le groupe  $\GL_{N_\rho}(\bbR)$ si l'on est dans le cadre des groupes classiques, 
et $\GL_{N_\rho}(\bbC)$ si l'on est dans le cadre des groupes unitaires.

Selon la forme intérieure $G$ et la dimension $N_\rho$ de $\rho$, 
 le groupe $G$ admet  ou pas un sous-groupe de Levi maximal standard $M$ isomorphe
à $\GL_{N_\rho} \times G'$, où $\mathbf G'$ est une forme intérieure de $\mathbf G^\flat$.
Par exemple, si $G=\U(p,q)$, la condition est que $\inf(p,q)\geq N_\rho$, et si c'est le cas, on a $G'=\U(p-N_\rho,q-N_\rho)$, et 
on a la même condition si $G=\SO(p,q)$. Si $\mathbf G$
est quasi-déployé, la condition est toujours vérifiée avec $\mathbf G'=\mathbf G^\flat$,  
et ceci fournit une injection
\begin{equation}\label{injGG}
\iota: \;    {}^LM=  \left(\widehat {\GL_{N_\rho}} \times \widehat{G^\flat}\right) \rtimes W_\bbR \hookrightarrow {}^LG
\end{equation}
de sorte que $\psi_G=\iota\circ \psi_M$ où $\psi_M: W_\bbR \rightarrow {}^LM$ est  construit à partir 
de  $\rho $ et $\psi_{G^\flat}$.
Ici $ \widehat {\GL_{N_\rho}}=\GL_N(\bbC)$ si $\GL_{N_\rho}=\GL_{N_\rho}(\bbR)$ et 
$ \widehat {\GL_{N_\rho}}=\GL_N(\bbC)\times \GL_N(\bbC)$ si $\GL_{N_\rho}=\GL_{N_\rho}(\bbC)$.

Comme le groupe $\mathbf G^\flat$ ne joue pas de rôle dans ce qui suit, on note plutôt  $\psi_{G'}$ pour $\psi_{G^\flat}$
du moins si la condition d'existence de $\mathbf G'$ est satisfaite. On vérifie facilement l'énoncé suivant.

\begin{rmq}\label{rmqA} Les groupes $A(\psi_G)$ et $A(\psi_{G'})$ sont naturellement isomorphes.        
\end{rmq}

Reprenons les notations de l'introduction, où pour un paramètre d'Arthur $\psi_G$ 
pour le groupe $G$, nous avons noté $\pi^A(\psi_G)$ la représentation unitaire de longueur finie
de $G\times A(\psi_G)$ attachée à $\psi_G$. On la note aussi $\pi^A(\psi_G,G)$.
On décompose maintenant cette représentation selon les caractères du groupe abélien fini $A(\psi)$ :
\begin{equation}\label{GA}
\pi^A(\psi,G)=\bigoplus_{\eta\in \widehat{A(\psi_G)}} \pi(\psi_G,\eta,G)\boxtimes \eta
\end{equation}
où les  $\pi(\psi_G,\eta,G)$ sont maintenant des représentations unitaires de longueur finie de $G$.

\begin{prop} \label{reducPaqArt} Avec les notations ci-dessus, si la condition d'existence de la forme intérieure
$\mathbf G'$ n'est pas vérifiée, on a $\pi^A(\psi_G)=0$. 

Si  la condition d'existence de la forme intérieure
$\mathbf G'$ est  vérifiée, 
soit $\eta \in \widehat{A(\psi_G)}$ et soient $\pi(\psi_G,\eta,G)$ et $\pi(\psi_{G'},\eta, G')$ les représentations semi-simples de $G$ et $G'$
respectivement attachées par Arthur ({\sl cf.} (\ref{GA}), où pour $\pi(\psi_{G'},\eta,G')$ on tient compte de la remarque 
ci-dessus). On a alors
\begin{equation}\label{IndpsiG}
\pi (\psi_G,\eta,G) =   \Ind_P^G  \left( \Pi^{\GL}_\rho\otimes \pi(\psi_{G'},\eta,G')\right),
  \end{equation}
où $P$ est un sous-groupe parabolique standard  maximal de $G$ de facteur de Levi $M$ isomorphe à $\GL_{N_\rho} \times G'$.
\end{prop}

Pour les groupes classiques, c'est  la proposition  4.3 de \cite{MR3}, énoncée sans démonstration.  

\dem 
Notons  $\pi^B(\psi_G,\eta,G)$ la représentation induite du membre de droite dans (\ref{IndpsiG}), et $\pi^B(\psi_G,G)=\oplus_{\eta\in \widehat {A(\psi_G)}}
\pi^B(\psi_G,\eta,G)\boxtimes \eta$. 
Nous avons  besoin de savoir que $\pi^B(\psi_G,G)$ est non nul si $G$ est quasi-déployé avant de pouvoir démontrer que cette représentation est $\pi^A(\psi_G,G)$.
Evidemment si $G$ est quasi-déployé, on a remarqué que  la condition d'existence de $\mathbf G'$ est toujours satisfaite, et que $\mathbf G'=\mathbf G^\flat$ 
est quasi-déployé. Ainsi   $\pi^B(\psi_G,G)$ est non nul si et seulement si $\pi^A(\psi_{G'}, G')$ est non nul. Or on  sait que
le paquet   $\Pi^A(\psi_{G'},G')$ est non nul, car il contient au moins les représentations dans le paquet de Langlands associé au paquet d'Arthur.
 Ceci montre l'assertion voulue. Considérons alors la représentation virtuelle stable $\pi^A(\psi_{G'},G')(s_{\psi'})=[\pi(\psi_{G'},G')]$ 
  ({\sl cf.} (\cite{MR3} (2.3.3)). Elle vérifie l'identité endoscopique    tordue (\cite{MR3} (3.2.4).
  Comme le transfert endoscopique tordu commute avec l'induction, on obtient que le transfert tordu de la représentation virtuelle
  $\pi^B(\psi_G,G)(s_\psi)$ est la trace tordue de l'induite  pour le parabolique standard de $\GL_N$ de Levi standard   $\GL_{N_\rho}\times \GL_{N'}$
  de la représentation $\Pi^\GL_\rho\boxtimes \Pi^{\GL}_{\psi'}$. Or cette induite est $\Pi^\GL_\psi$, d'après la définition des paquets d'Arthur des groupes linéaires
  (voir \cite{AMR}, section 3.2).
  D'autre part   $\pi^A(\psi_G,G)(s_\psi)=[\pi(\psi_G,G)]$  vérifie aussi  l'identité endoscopique tordue.
On obtient donc que $\pi^A(\psi_G,G)(s_\psi)=\pi^B(\psi_G,G)(s_\psi)$, puisque ces deux représentations virtuelles stables ont même transfert endoscopique tordu.

Nous allons démontrer que  $\pi^B(\psi_G,G)$  vérifie aussi les identités endoscopiques ordinaires (\cite{MR3} (2.3.5)). On ne suppose plus
que $G$ est quasi-déployé, mais l'on suppose que la condition d'existence de $\mathbf G'$ est satisfaite, car sinon, il est clair que
$\pi^A(\psi_G,G)=0$.
Soit  ${\bf H}=(H,x,\xi: {}^LH\rightarrow {}^LG, \ldots)$ une  donnée endoscopique elliptique de $G$
({\sl cf.} \cite{Art13}) telle que $\psi_G$ se factorise 
par le groupe dual de $H$ et on fixe une telle factorisation $\psi_G=\xi \circ\psi_x$. 
En particulier  l'élément $x\in {}^LG$
 s'identifie à un élément du commutant de $\psi_G$. Il faut alors démontrer qu'il 
  existe une donnée endoscopique elliptique ${\bf H}'=(H',x',\ldots)$ de $G'$,  tel que l'élément $x'$ de 
 cette donnée soit dans le centralisateur de $\psi_{G'}$ et tel que le transfert de la distribution stable 
 associée à ${\bf H}$ et à la factorisation de $\psi_G$ soit l'induite du produit tensoriel 
  des données analogues pour $\psi_{G'}$ et   ${\bf H'}$ et de  la représentation $ \Pi^{GL}_\rho$. 
  Expliquons maintenant comment construire explicitement cette donnée endoscopique ${\bf H}'$.
 Comme il est loisible de le faire ici, on suppose que $x$ vérifie $x^2=1$. 
 On décompose alors $\psi$ en $\psi_{+}\oplus \psi_-$ suivant les valeurs propres de $x$. 
 On remarque que l'on a aussi une décomposition analogue pour $\psi'$ et pour $\rho$. On a alors 
$$\psi_{+}= \rho_{+}\oplus \rho_{+}^* \oplus \psi'_{+}$$ 
et une décomposition analogue avec $+$ remplacé par $-$.
C'est ici qu'a servi l'hypothèse sur la mauvaise parité des composantes de $\rho$, pour que le dual 
de $\rho_+$ apparaisse lui aussi dans  dans l'espace propre de valeur propre $+1$.
Notons $N_{\rho_\pm}$  les dimensions des représentations $\rho_\pm$, 
et $\Pi^\GL_{\rho_\pm}$ la représentation de $\GL_{N_{\rho_\pm}}$ associée à ce paramètre.
On a bien s\^ur  $N_\rho=N_{\rho_+}+N_{\rho_-}$ et $\Pi^\GL_{\rho}$ est l'induite parabolique  
de $\Pi^\GL_{\rho_+}\boxtimes \Pi^\GL_{\rho_-}$.

Ainsi il existe un sous-groupe de Levi 
$$M_H \simeq \left( \GL_{N_{\rho_+}} \times M^+\right)\times \left(\GL_{N_{\rho_-}}  \times M^-\right)$$ 
de  $H$ tel que $\psi_x$ se factorise par le $L$-groupe de  $M_H$ 
et la représentation virtuelle  stable  $[\pi(\psi_x,H)] $ de $H$ associée à $\psi_x$  
 est une induite à partir de ce Levi. Notons $H'$ le facteur    $M^+\times M^-$
 de $M$ : c'est un groupe endoscopique pour $G'$, s'inscrivant dans une donnée endoscopique 
 ${\bf H}'=(H',x', \xi',\ldots)$ de $G'$ et le paramètre d'Arthur $\psi_{G'}$ se factorise
 en $\xi'\circ \psi_{x'}$. L'élément $x'$ est dans le centralisateur de $\psi_{G'}$, on
 peut le prendre tel que ${x'}^2=1$ et $\psi'={\psi'}^+\oplus {\psi'}^-$ est la décomposition de 
 $\psi'$ selon les valeurs propres $\pm 1$ de $x'$. 
 Partons de la représentation stable $[\pi(\psi_{x'},H')]$ 
associée à $\psi_{x'}$. On peut d'abord considérer son transfert endoscopique 
vers $G'$, puis induire vers $G$ avec $\Pi^{\GL}_{\rho}$ : 
 $$  \Ind_{P=MN}^G \left( \Pi^\GL_{\rho}\boxtimes  \mathrm{Trans}_{H'}^{G'} ( [\pi(\psi_{x'},H')] )\right)   $$
où  $  \mathrm{Trans}_{H'}^{G'}$ désigne le transfert endoscopique (spectral) du groupe endoscopique $H'$ de $G'$ vers $G'$.
Or, ceci est égal à 
$$  \Ind_{P=MN}^G \left( \Pi^\GL_{\rho}\boxtimes  \pi^A(\psi_{G'},H')(s_{\psi'}x')  )\right)  = \pi^B(\psi_G,G)(s_\psi x). $$
Le fait que le transfert commute à l'induction nous dit que l'on obtient le même résultat 
en prenant le produit tensoriel extérieur avec  $\Pi^\GL_{\rho_+}$ et $\Pi^\GL_{\rho_-}$ pour  obtenir une représentation virtuelle de $M_H$ que l'on induit    
vers $H$, puis en prenant ensuite le transfert 
endoscopique de $H$ vers $G$ : 
$$ \mathrm{Trans}_{H}^{G}( \Ind_{P_h=M_HN_H}^H \left(   \Pi^\GL_{\rho_+}\boxtimes   \Pi ^\GL_{\rho_-} \boxtimes   [\pi(\psi_{x'},G')] \right).$$
Or ceci est égal à 
$$ \mathrm{Trans}_{H}^{G}(   [\pi(\psi_{x},H)] ) = \pi^A(\psi_G,G)(s_\psi x) .$$
On obtient donc que $ \pi^B(\psi_G,G)(s_\psi x) =  \pi^A(\psi_G,G)(s_\psi x) $. 
Comme on a ceci pour tout $x\in A(\psi_G)$, que $\mathbf G$ soit ou non quasi-déployé, on en déduit par inversion de Fourier
que  $ \pi^B(\psi_G,G)=  \pi^A(\psi_G,G)$. 
\qed

\subsection{Irréductibilité de l'induite parabolique pour les groupes unitaires} \label{irredU}

Nous allons reformuler la proposition \ref{reducPaqArt}  de manière un peu plus explicite pour les groupes unitaires, en y ajoutant un résultat 
d'irréductibilité des induites paraboliques.

\begin{thm}\label{ThmReduc}
 Soient $\psi_G$ un $A$-paramètre pour $G=\U(p,q)$,  $\psi$ sa restriction à $\bbC^\times \times \SL_2(\bbC)$
 comme en (\ref{decApar}) et  $\psi_{bp}$ la partie de bonne parité de ce paramètre. Soient $N=p+q$ et  $N_{bp}$ 
comme en (\ref{nbp}). En particulier $N-N_{bp}$ est pair, et l'on pose $a_{mp}= \frac{N-N_{bp}}{2}$, c'est le cardinal de l'ensemble $\caE'(\psi)$. On a alors :

(i) si $\inf(p,q)< a_{mp}$,  alors $\pi^A(\psi_G, G)=0$.

(ii) Si $\inf(p,q)\geq a_{mp}$  on pose $p_{bp}= p- a_{mp}$, $q_{bp}=q- a_{mp}$ et  $\psi_{bp}$ est la restriction  à $\bbC^\times \times \SL_2(\bbC)$
d'un paramètre  $\psi_{G_{bp}}$ pour $G_{bp}=\U(p_{bp},q_{bp})$. 
On a donc une représentation 
unitaire   $\pi^A(\psi_{G_{bp}},G_{bp})$ de $\U(p_{bp},q_{bp})\times A(\psi_{G_{bp}})$, qui s'écrit 
\[   {\pi}^A(\psi_{G_{bp}}, G_{bp})=\bigoplus_{\eta\in \widehat{A(\psi_{G_{bp}})}}  \pi(\psi_{G_{bp}},\eta,G_{bp}) \boxtimes \eta. \] 
Alors $\pi^A(\psi_G, G  )$ s'écrit 
$  {\pi}^A(\psi_{G}, G)=\bigoplus_{\eta\in \widehat{A(\psi)}}  \pi(\psi_G,\eta,G) \boxtimes \eta$,  
avec pour tout $\eta\in A(\psi_G)$ (rappelons qu'en vertu de la remarque \ref{rmqA}, on peut identifier  $A(\psi_{G_{bp}})$ et  $A(\psi_{G})$), 
\[  \pi(\psi_G, \eta, G) = \Ind_P^G \left(   \left( \boxtimes_{(t,s,a)\in \mathcal{E}'(\psi)}\chi_{t,s,a}\right)\boxtimes \pi(\psi_{G_{bp}},\eta, G_{bp}) \right)\]
pour le sous-groupe parabolique standard $P$ de $\U(p,q)$ dont le sous-groupe de  Levi est 
$\times_{(t,s,a)\in \mathcal{E}' (\psi)} \GL(a,\mathbb{C}) \times \U(p_{bp},q_{bp}) $. 

De plus, si $\tau$ est une sous-représentation irréductible de $\pi(\psi_{G_{bp}},\eta, G_{bp}) $, alors 
$\Ind_P^G \left( \left(  \boxtimes_{(t,s,a)\in \mathcal{E}'(\psi)}\chi_{t,s,a}\right)\boxtimes \tau \right)$
 est irréductible.

\end{thm}

\dem
Seule la dernière assertion est nouvelle par rapport à la proposition.
Les représentations $\tau$ de la dernière assertion du théorème sont les 
représentations $\scrA_{\underline{d}}(\psi_{bp})$ de la section \ref{bonneparite}  attachée à la partie de bonne parité du paramètre.
Il est démontré dans \cite{Mat96} 3.2.2  que l'induite parabolique de
\[ 
\left(\boxtimes_{(t,s,a)\in \mathcal{E}'(\psi)}\chi_{(t,s,a)} \right)\boxtimes \scrA_{\underline{d}}(\psi_{bp})  \]
est irréductible, ce qui est exactement l'assertion voulue. 
\qed

\subsection{Irréductibilité de l'induction parabolique pour les groupes classiques}

Nous complétons maintenant la proposition \ref{reducPaqArt}  pour les groupes classiques
en y ajoutant le fait que, comme pour les groupes unitaires,  l'induction parabolique préserve l'irréductibilité.
Ce résultat avait été énoncé sans démonstration dans  \cite{MR3},  Théorème 4.4.

\begin{thm}\label{redCla}
On se place dans les hypothèses de la proposition \ref{reducPaqArt}. On suppose que la condition d'existence de la forme intérieure
$\mathbf G'$ est  vérifiée.  Soit $\eta \in \widehat{A(\psi_G)}$. 
 Si $\tau$ est une sous-représentation irréductible de $\pi(\psi_{G'},\eta, G') $, alors 
$ \Ind_P^G  \left( \Pi^{\GL}_\rho\boxtimes \tau \right) $
 est irréductible.
\end{thm}

Rappelons que pour un groupe classique $G$, la mauvaise parité est : impaire si $\widehat G$ est un groupe spécial orthogonal, 
paire si  $\widehat G$ est un groupe symplectique.
On écrit la partie de mauvaise parité $\rho$   ({\sl cf. \cite{MR3}, \S 4.1})  sous la forme 
\[\rho=\bigoplus_{i=1,\ldots,\ell} \delta_{t_i,s_i}\boxtimes R[a_i] \oplus  \bigoplus_{j=1,\ldots,\ell'} \eta_{\epsilon_j,s_j}\boxtimes R[a'_j].\] 
Dans la première somme, $t_i\in \bbZ_{>0}$, $s_i\in i\bbR$,  et $\delta_{t_i,s_i}$ est le paramètre de Langlands de la série discrète de caractère infinitésimal 
$\left(\frac{t_i+s_i}{2},\frac{-t_i+s_i}{2} \right)$ de $\GL_2(\bbR)$, et si $s_i=0$, alors $t_i+a_i-1$ est de mauvaise parité. Dans la seconde somme
$\epsilon_j\in \{\pm 1\}$, $s_j\in i\bbR$, et $ \eta_{\epsilon_j,s_j}$ est le paramètre de Langlands du caractère $x\mapsto \sgn(x)^{\frac{1-\epsilon_j}{2}}\vert x\vert^{s_j}$ 
 de $\GL_1(\bbR)$, et si  $s_j=0$, alors $a'_j-1$ est de mauvaise parité.
  On note encore $\delta_{t_i,s_i}$  et $ \eta_{\epsilon_j,s_j}$ les représentations de $\GL_2(\bbR)$ et $\GL_1(\bbR)$ dont ce sont les paramètres.
Pour chaque indice $i$, on considère la représentation de Speh, notée  $\mathbf{Speh}(\delta_{t_i,s_i},a_i  )$ de $\GL_{2a_i}(\bbR)$
 qui est irréductible et unitaire, et pour chaque
indice $j$, le caractère unitaire  $\eta_{\epsilon_j,s_j} \circ \det$ de $\GL_{a'_j}(\bbR)$. 
La représentation  $\Pi^{\GL}_\rho$ est alors obtenue par induction parabolique irréductible à partir de la représentation
\[  \left( \boxtimes_{i=1,\ldots,\ell} \,  \mathbf{Speh}(\delta_{t_i,s_i},a_i  ) \right) \boxtimes \left( \boxtimes_{j=1,\ldots,\ell}\,  \eta_{\epsilon_j,s_j} \circ \det \right) \]
du facteur de Levi 
\[    \left( \times_{i=1,\ldots,\ell} \, \GL_{2a_i}(\bbR)   \right) \times \left( \times_{j=1,\ldots,\ell}\,  \GL_{a'_j}(\bbR)\right) \]      
de $\GL_{N_\rho}(\bbR)$.
\, 

Soit $G_0$ le groupe de  même type que $G$ tel que  $\GL_{N_\rho}(\bbR) \times G_0$ est un   sous-groupe de Levi d'un parabolique $P$ de $G$.
 On note $\tau_0$ une représentation unitaire irréductible de $G_0$. On note $N_0$ son rang.
 On suppose que le caractère infinitésimal de $\tau_0$ 
 a bonne parité, c'est-à-dire qu'il  est formé d'entiers si la demi-somme des racines positive de $G_0$ est formée d'entiers et est formé de demi-entiers non entiers sinon.
  Le fait que $\tau_0$ soit unitaire n'est absolument pas nécessaire mais simplifie légèrement la preuve.
  Le théorème résulte alors de la proposition plus générale suivante.
 \begin{prop} \label{redCla2}  La représentation induite
$\Ind _P^G  \left( \Pi^{GL}_\rho \boxtimes \tau_0\right) $  est irréductible.
 \end{prop}
 
 \dem  On note $\pi=\Ind _P^G  \left( \Pi^{GL}_\rho \boxtimes \tau_0\right) $. Nous allons scinder la démonstration en plusieurs étapes.

 {\sl Première étape}.
 Pour tout $j\in [1,\ell']$, on remplace $\eta_{\epsilon_j,s_j}\boxtimes R[a'_j]$ dans la partie de mauvaise parité $\rho$ par 
 \[ (\eta_{\epsilon_j,s_j}\boxtimes R[a'_j)]\oplus (\eta_{-\epsilon_j,s_j}\boxtimes R[a'_j]) . \]
 On obtient ainsi un paramètre $\rho^\sharp$ de dimension plus grande, et toujours de mauvaise parité, et il est clair 
 que la proposition est vraie pour $\rho$ si elle l'est pour $\rho^\sharp$.
Remarquons que l'on peut poser  $ \eta_{\epsilon_j,s_j}\oplus \eta_{-\epsilon_j,s_j} =\delta_{0,s_j}$ (la limite de séries discrètes $\delta_{0,s_j}$de $\GL_2(\bbR)$
est l'induite de  $\eta_{\epsilon_j,s_j}\boxtimes  \eta_{-\epsilon_j,s_j}$).
Ainsi, on peut supposer que $\rho=\bigoplus_{i=1,\ldots,\ell} \delta_{t_i,s_i}\boxtimes R[a_i]  $, mais il est maintenant possible que certains $t_i$ soient nuls.
La dimension $N_\rho$ de $\rho$ est  paire, et l'on pose $N'_\rho=N_\rho/2=\sum_i a_i$.
 On change maintenant légèrement les notations, $P=MN$ désigne maintenant le sous-groupe parabolique 
standard de $G$ de facteur de Levi isomorphe à  $\left( \times_i \GL_{2a_i}(\bbR) \right)\times G_0$.
La représentation $\pi$ est donc avec ces notations
 \[\pi =\Ind _P^G  \left(\left(  \boxtimes_i\mathbf{Speh}\, (\delta_{t_i,s_i},a_i)   \right) \boxtimes \tau_0\right). \]
 \
 
  {\sl Deuxième  étape}.
Les représentations  $\mathbf{Speh}\, (\delta_{t_i,s_i},a_i) $ sont obtenues à partir du 
caractère $\chi_{t_i,s_i,a_i}$ de $\GL_{a_i}(\bbC)$ par induction cohomologique. Ceci est bien connu, voir par exemple \cite{KnVo}, p. 586, et l'on uitilise
ici la version normalisée de l'induction cohomologique {\sl loc. cit.} (11.150b).
Ainsi, $\pi$ est obtenue en deux étapes, d'abord une induction cohomologique, puis une induction parabolique.
Le théorème de transfert du chapitre 11 de {\sl loc. cit.} permet d'échanger l'ordre de ces deux inductions. Une référence commode pour cela est 
\cite{Mat04}, Thm. 2.2.3 qui nous donne exactement l'énoncé voulu. Donnons les détails.
Matumoto introduit la terminologie de $\sigma\theta$-paire pour  un couple $(\frp,\frqqq)$ de sous-algèbres paraboliques
de $G$. La sous-algèbre $\frp$ est la complexifiée de l'algèbre de Lie d'un sous-groupe parabolique 
$P$ de $G$, qui ici a été fixée à la fin de la première étape.
La sous-algèbre $\frqqq$ est une sous-algèbre parabolique $\theta$-stable. On pose $L=\mathrm{Norm}_G(\frqqq)$.
On choisit ici cette sous-algèbre de sorte que d'une part $L$ soit  isomorphe à $U(N'_\rho,N'_\rho)\times G_0$ et que d'autre part
la condition S2 de la définition 2.2.1 de \cite{Mat04} soit vérifiée, c'est-à-dire que $\frp\cap \frqqq$ contient une sous-algèbre de Cartan $\sigma$ et $\theta$-stable
de $\frg$. Il est facile de vérifier que l'on trouve une telle sous-algèbre $\frqqq$, en partant d'un sous-groupe de Cartan isomorphe
à $(\bbC^\times)^{N'_\rho}\times \U(1)^{N_0}$ de $M$. On adopte les notations de Matumoto, qui sont usuelles ($\frp=\frmm\oplus \frn$, $\frqqq=\frl\oplus \fru$, etc).
On a ainsi 
\[  \pi=\Ind _P^G  \left(   \left( {}^n\caR_{\frqqq \cap \frmm,L\cap M \cap K}^{\frmm,M\cap K}   \right)^{\dim(\fru\cap \frmm\cap \frk)}   \left(    \left(\boxtimes_i
\chi_{t_i,s_i,a_i}\right) \boxtimes \tau_0\right)  \right) . \]
Remarquons que ici $L\cap M$ est isomorphe à $(\times_i \GL_{a_i}(\bbC))\times G_0$.
Le résultat de Matumoto nous permet  alors d'écrire, sous certaines conditions de positivité des paramètres,
\[ \pi=  \left( \caR_{\frqqq ,L\cap  K}^{\frg,K}   \right)^{\dim(\fru \cap \frk)}  \left(   \Ind_ {P\cap L}^L\left(  \left( \boxtimes \chi_{t_i,s_i,a_i}\right) \boxtimes \tau_0) 
\otimes \bbC_{-\delta(\fru)} 
  \right)  \right) . \]
Ici, on a simplifié la formulation de Matumoto, en utilisant le fait que le groupe $L\cap M$ est connexe, et ainsi
le caractère noté $\bbC_{\delta(\bar \frn\cap \fru)'}$ par Matumoto venant de  la normalisation subtile  des foncteurs d'induction cohomologique
dans \cite{KnVo} coïncide ici avec  $\bbC_{\delta(\bar \frn\cap \fru)}$ (en général ces deux caractères sont seulement égaux sur la composante neutre de $L\cap M$).
D'autre part, si l'on suppose les $t_i$ suffisamment grand, les hypothèses de positivité dans le théorème de Matumoto sont vérifiées.

On veut montrer que $\pi$ est irréductible. Or, avec les $t_i$ suffisamment grand,
l'induction cohomologique  $\left( \caR_{\frqqq ,L\cap  K}^{\frg,K}   \right)^{\dim(\fru \cap \frk)} $ se fait dans le good range et préserve donc l'irréductibilité.
Il suffit alors de démontrer que   $\Ind_ {P\cap L}^L\left(  \left( \boxtimes \chi_{t_i,s_i,a_i}\right) \boxtimes \tau_0) 
\otimes \bbC_{-\delta(\fru)}   \right) $ est irréductible. Or 
\[\Ind_ {P\cap L}^L\left(  \left( \boxtimes \chi_{t_i,s_i,a_i}\right) \boxtimes \tau_0) 
\otimes \bbC_{-\delta(\fru)}   \right) =\left( \Ind_{P'}^{\U(N'_\rho,N'_\rho)}   \left( \left(\boxtimes \chi_{t_i,s_i,a_i} \right)\otimes  \bbC_{-\delta(\fru)}   \right)  
 \right) \boxtimes \tau_0 \]
Ici $P'$ est un sous-groupe parabolique de $\U(N'_\rho,N'_\rho)$ de facteur de Levi isomorphe à $\times_i\GL_{a_i}(\bbC)$ (ce dernier se plonge 
naturellement dans $\GL_{N_\rho}(\bbC)$ et l'on voit  $\GL_{N_\rho}(\bbC)$ comme le Levi du  parabolique de Siegel). 
On est ramené à montrer l'irréductibilité de $\Ind_{P'}^{\U(N'_\rho,N'_\rho)}   \left( \left(\boxtimes \chi_{t_i,s_i,a_i} \right)\otimes  \bbC_{-\delta(\fru)}   \right) $. 
Ceci découle de la dernière assertion du théorème \ref{ThmReduc} car 
 l'hypothèse de mauvaise parité pour le groupe unitaire $\U(N'_\rho,N'_\rho)$ est vérifiée, comme on le montre ci-dessous.

On calcule le caractère de torsion $\bbC_{-\delta(\fru)}$.
Posons $\epsilon_G=0$ si $G$ est un groupe orthogonal pair, 
$\epsilon_G=1$ si $G$ est un groupe  symplectique, 
et $\epsilon_G=\frac{1}{2}$ si $G$ est un groupe orthogonal impair.
Dans le système de coordonnées usuelles pour $G$, on a
\[\delta(\fru)= (\underbrace{ N-N'_\rho-\frac{1}{2}+\epsilon_G,\ldots, N-N'_\rho-\frac{1}{2}+\epsilon_G}_{N_\rho},\underbrace{0,\ldots ,0}_{N_0}).  \]

L'hypothèse de mauvaise parité des $(t_i,s_i,a_i)$ pour $G$ est : soit $s_i\neq 0$, soit $s_i=0$ et 
$\frac{t_i+a_i-1}{2}+\epsilon_G\in \frac{1}{2}\bbZ\setminus \bbZ$.
Dans les deux cas,  $(t_i+N-N'_\rho-\frac{1}{2}+\epsilon_G,s_i,a_i)$ est de mauvaise parité pour $\U(N'_\rho,N'_\rho)$, car dans le second cas
\[ \frac{t_i+a_i-1}{2}  +  N-N'_\rho-\frac{1}{2}+\epsilon_G+ \frac{2N'_\rho-1}{2}  \in \frac{1}{2}\bbZ\setminus \bbZ.\]

Ceci termine la deuxième étape où l'on a établi le résultat voulu si la condition que les $t_i$ soient suffisamment grands.

\,

{\sl Troisième étape}. On ne suppose plus ici que les $t_i$ sont suffisamment grand, mais on choisit 
$T$ suffisamment grand pour que la représentation $\pi_T$ obtenue en remplaçant les 
$t_i$ par $t_i+T$ soit irréductible d'après la deuxième étape. C'est-à-dire que l'on pose
 \[\pi_T =\Ind _P^G  \left(\left(  \boxtimes_i\mathbf{Speh}\, (\delta_{t_i+T,s_i},a_i)   \right) \boxtimes \tau_0\right) \]
\[ =\Ind _P^G  \left(   \left( {}^n\caR_{\frqqq \cap \frmm,L\cap M \cap K}^{\frmm,M\cap K}   \right)^{\dim(\fru\cap \frmm\cap \frk)}   \left(    \left(\boxtimes_i
\chi_{t_i+T,s_i,a_i}\right) \boxtimes \tau_0\right)  \right) 
 \]

On veut montrer que $\pi$ s'obtient de $\pi_T$ par un foncteur de translation (\cite{KnVo}, Chapter 7), c'est-à-dire
que l'on obtient $\pi$ en tensorisant $\pi_T$ par une représentation de dimension finie $F$ de $G$, et l'on projette
sur la composante de caractère infinitésimal généralisé le caractère infinitésimal de $\pi$.
 On note $n_0$ le rang de $G_0$. Le rang de $G$ est donc $n_0+N_\rho=n_0+2N'_\rho$.

On choisit une sous-algèbre de Cartan $\frh$ de $\frg$ et un système de racines positives $\Delta^+(\frg,\frh)$ de $\frh$ dans $\frg$.
On considère la représentation de dimension finie $F$, de $G$ de plus haut poids 
 \[  (\underbrace{T/2,\ldots,T/2}_{2N'_\rho} , \underbrace{0,\ldots,0}_{n_0}). \]
 On remarque que cette représentation admet le poids extrémal
 \[
 \mu_0:=(\underbrace{-T/2, \ldots,- T/2}_{N'_\rho}, \underbrace{T/2, \ldots, T/2}_{N'_\rho},  \underbrace{0,\ldots,0}_{n_0}). \]

 On note $\mathfrak{p}'=\frmm'\oplus\frn'$ la sous-algèbre parabolique  de $\frg$ qui stabilise le sous-espace poids $\mu_0$  de $F$  
 (qui est de dimension $1$ car $\mu_0$ est extrémal). Son  facteur de Levi   $\frmm'$ est  isomorphe à  
 $\frg\frl_{N'_\rho}(\bbC)\times \frg \frl_{N'_\rho}(\bbC) \times \frg_0$. On note $\mathfrak{q}'=\frl'\oplus \fru'$ la sous-algèbre parabolique de $\frg$ contenue  dans $\mathfrak{p}'$ 
 dont la sous-algèbre de Levi  $\frl'$ est isomorphe à 
 \[\frg\frl_{a_1}(\mathbb{C})\times\cdots\times \frg\frl_{a_\ell}(\mathbb{C})\times \frg\frl_{a_\ell} (\mathbb{C}) \times \cdots \times 
 \frg\frl_{a_1}(\mathbb{C})\times \frg_0.
 \]
 Ceci ne détermine pas  pas $\mathfrak{q}'$ : on a $\frl'\subset \frmm'$, $\frn'\subset \fru'$,  $\fru'=\frn'\oplus (\frmm'\cap \fru')$ et
l'on fixe $ \frmm'\cap \fru'$ en demandant que pour toute racine $\beta\in \Delta(\fru',\frh)$, 
\begin{equation}\label{condpos} \beta( \underbrace{t_1,\ldots,t_1}_{a_1}, \ldots, \underbrace{t_\ell,\ldots,t_\ell}_{a_\ell} \underbrace{-t_\ell,\ldots,t_\ell}_{a_\ell},\ldots,   \underbrace{-t_1,\ldots,-t_1}_{a_1} )
\leq 0.\end{equation}

 On note $\tilde{\lambda}$ le caractère infinitésimal de $\pi$ et $\tilde{\lambda}_T$ celui de $\pi_T$ que l'on voit aussi comme éléments de $\frh^*$. 
On veut montrer que $\pi$ s'obtient de $\pi_T$ par le foncteur de translation défini par $F$, c'est-à-dire que 
\[  \pi=(\pi_T\otimes F)_{\tilde \lambda} \]
où $(\, . \,)_{\tilde \lambda}$ dénote la projection d'un module $\frZ(\frg)$-fini sur sa composante primaire pour le  caractère infinitésimal 
 généralisé défini par  $\tilde \lambda$.

  Pour  cela, on suit  \cite{VogDs} Prop. 4.7 en modifiant convenablement les hypothèses. 
 On commence par montrer un analogue du lemme 4.8  de Vogan.  
 \begin{lemme}Soit $\mu\in \frh^*$ le plus haut poids d'une composante irréductible de la restriction de $F$ à $\frl'$. On suppose que pour un certain 
   $w$ dans le groupe de Weyl de $G$, on a  
 \[
 (\tilde{  \lambda}_T+\mu)=w.\tilde{ \lambda},
 \]
Alors $\mu=\mu_0$. 
 \end{lemme}

\dem 
 On pose  $\frh_\frl=\frh\cap [\frl',\frl']$ et l'on note $\frz_\frl$ le centre de $\frl'$. On a alors $\frh=\frh_\frl\oplus \frz_\frl$ et $\frh^*= \frh_\frl^*\oplus \frz_\frl^*$. 
 C'est une somme directe orthogonale.  Si $\nu\in \frh^*$,  on  note $\nu=\nu_\frl \oplus \nu_\frz$ sa décomposition selon cette somme orthogonale. 
 Le poids $\mu$ s'écrit $\mu_0+\sum_{\beta\in \Delta(\frn',\frh)} m_\beta \, \beta$, avec les coefficients $m_\beta$ entiers négatifs.
 On a donc 
 \[ (\tilde \lambda_T +\mu)_\frz= (\tilde \lambda_T +\mu_0)_\frz+ \sum_{\beta\in \Delta(\frn',\frh)} m_\beta \, \beta_{\vert \frz_\frl} \]

 Dans le système de coordonnées choisi, $\tilde \lambda_T$ s'écrit 
 {\tiny\[  \frac{1}{2}( \underbrace{t_1+T+(a_1-1), \ldots , t_1+T-(a_1-1)}_{a_1} , \ldots,  \underbrace{t_\ell+T+(a_\ell-1), \ldots , t_\ell+T+(a_\ell-1)}_{a_\ell} , \]
\[ \underbrace{-t_\ell-T+(a_\ell-1), \ldots , -t_\ell-T-(a_\ell-1)}_{a_\ell}, \underbrace{-t_1-T+(a_1-1), \ldots , -t_1-T-(a_1-1)}_{a_1} , *, \ldots ,*   )  \]}
 où sur les dernière coordonnées, il apparaît le caractère infinitésimal de $\tau_0$ que nous n'explicitons pas.
 Ainsi  $\tilde \lambda_T+\mu_0$ s'écrit 
  {\tiny\[  \frac{1}{2}( \underbrace{t_1+(a_1-1), \ldots , t_1-(a_1-1)}_{a_1} , \ldots,  \underbrace{t_\ell+(a_\ell-1), \ldots , t_\ell+(a_\ell-1)}_{a_\ell} , \]
\[ \underbrace{-t_\ell+(a_\ell-1), \ldots , -t_\ell-(a_\ell-1)}_{a_\ell}, \underbrace{-t_1-+(a_1-1), \ldots , -t_1-(a_1-1)}_{a_1} , *, \ldots ,*   )  \]}
 et $\tilde \lambda_T+\mu_0=\tilde \lambda$, d'où 
  \[ (\tilde \lambda_T +\mu)_\frz= \tilde \lambda_\frz+ \sum_{\beta\in \Delta(\frn',\frh)} m_\beta \, \beta_{\vert \frz_\frl} \]
et en utilisant (\ref{condpos}), on en conclut
\begin{equation}\label{iNeg}  \Vert  (\tilde \lambda_T +\mu)_\frz\Vert \leq    \Vert   \tilde \lambda_\frz  \Vert \end{equation}
 avec égalité si tous les $m_\beta$ sont nuls.
 
D'autre part $\tilde{\lambda}_T+\mu$ et $\tilde{\lambda}=\tilde \lambda_T+\mu_0$ ont même projection sur $\frh_\frl$ car $\mu-\mu_0\in \frh^*_\frz$.   
Ainsi  $  \Vert  (\tilde \lambda_T +\mu)_\frl \Vert =    \Vert   \tilde \lambda_\frl  \Vert $.
Or l'hypothèse  force $\tilde{\lambda}_T+\mu$ et $\tilde{\lambda}$ à avoir même norme, et il y a donc égalité dans (\ref{iNeg}), d'où $\mu=\mu_0$. \qed 

\
 
 Pour en déduire le fait que la translation de $\pi_T$ est bien $\pi$, on raisonne comme dans \cite{VogDs}, où la proposition 4.7 est déduite du lemme 4.8.
 Ici, on utilise le fait que pour toute représentation $\tau$ de $M$ et pour toute représentation de dimension finie $F$ de $G$, 
 $\Ind_{P}^G(\tau)\otimes F=\Ind_P^G(\tau\otimes F_{\vert M})$. On applique ceci à 
 $\tau = \left( {}^n\caR_{\frqqq \cap \frmm,L\cap M \cap K}^{\frmm,M\cap K}   \right)^{\dim(\fru\cap \frmm\cap \frk)}   \left(    \left(\boxtimes_i
\chi_{t_i+T,s_i,a_i}\right) \boxtimes \tau_0\right) $ et à $F$ comme ci-dessus, et l'on obtient 
\[\pi_T\otimes F =\Ind _P^G  \left(   \left( {}^n\caR_{\frqqq \cap \frmm,L\cap M \cap K}^{\frmm,M\cap K}   \right)^{\dim(\fru\cap \frmm\cap \frk)}   \left(    \left(\boxtimes_i
\chi_{t_i+T,s_i,a_i}\right) \boxtimes \tau_0\right)  \otimes F_{\vert M}  \right). \]

 Ensuite, on utilise le fait que l'induction cohomologique a lieu dans le weakly fair range, où il y a annulation 
 des foncteurs d'induction cohomologique en degré différent de $\dim(\fru\cap \frmm\cap \frk)$, ce qui nous permet de 
 remplacer le foncteur $\left( {}^n\caR_{\frqqq \cap \frmm,L\cap M \cap K}^{\frmm,M\cap K}   \right)^{\dim(\fru\cap \frmm\cap \frk)}$
 par $\caR:=\sum_i (-1)^i\, \left( {}^n\caR_{\frqqq \cap \frmm,L\cap M \cap K}^{\frmm,M\cap K}   \right)^{\dim(\fru\cap \frmm\cap \frk)-i}$
 (on obtient alors une égalité de représentations virtuelles). Ceci nous permet d'utiliser \cite{Vgreen}, Lemma 7.2.9 (b), et l'on a alors 
 \[\pi_T\otimes F =\Ind _P^G  \left( \caR   \left(    \left(\boxtimes_i
\chi_{t_i+T,s_i,a_i}\right) \boxtimes \tau_0\right)  \otimes F_{\vert M}  \right)= \Ind _P^G  \left( \caR  \left( \left(    \left(\boxtimes_i
\chi_{t_i+T,s_i,a_i}\right) \boxtimes \tau_0\right)  \otimes F_{\vert M\cap L}  \right)\right) . \]

 On remarque alors, avec les notations employées, que l'on a en fait $\frl'=\frl\cap \frmm$  et $\frqqq'=\frl'\oplus \frn$.
 On conclut en remarquant que les contributions à la projection de $\pi_T\otimes F$
 sur la composante primaire $\tilde \lambda$ proviennent des composantes de la restriction de $F$ à $\frl'$ 
de  plus haut poids $\mu$ vérifiant l'hypothèse du lemme, et que ceci donne alors le résultat voulu.

 \ 
 
{\sl Quatrième étape}. 
  La translation préserve l'irréductibilité si l'image de l'algèbre enveloppante dans l'algèbre des endomorphismes $G$-finis de la représentation 
  induite du caractère de  $\mathfrak{p}$ est surjective. Un critère pour cela est que l'orbite de Richardson du parabolique $P$ 
  soit de fermeture normale et que l'application moment soit birationnelle. Les orbites décrites par Barbasch en \cite{B1} 14.5 
  vérifient ces critères. La description de Barbasch en termes d'induite de Richardson s'applique directement pour nous: dans le cas des groupes symplectiques,
   il suffit  que $G_0$ soit trivial, et pour les groupes orthogonaux à l'inverse, il suffit que le rang de $G_0$ soit grand par rapport aux $a_i$, 
  et que la représentation $\tau_0$ soit de dimension finie. Comme on suppose $\tau_0$ unitaire, si $\tau_0$ n'est pas de dimension
  $1$, cela veut dire que $G_0$ est un groupe compact, mais nous n'allons utiliser que le cas où $\tau_0$ est la représentation triviale.
 Ainsi, on obtient  l'irréductibilité de $\pi$  à condition que $\tau_0$  soit la représentation triviale du groupe trivial si $G$ est symplectique, et la représentation 
triviale  d'un groupe $G_0$ de rang suffisamment grand si $G$ est orthogonal.

 \

 {\sl Cinquième étape}.
 On revient à l'énoncé général de la proposition en utilisant toutefois la réduction effectuée dans la première étape. On a donc 
\[\pi =\Ind _P^G  \left(\left(  \boxtimes_i\mathbf{Speh}\, (\delta_{t_i,s_i},a_i)   \right) \boxtimes \tau_0\right). \]
  
  L'idée est de construire un opérateur d'entrelacement de $\pi$ dans dans un module standard en position de Langlands négative,
 l'irréductibilité de $\pi$ étant alors conséquence de l'injectivité de  cet opérateur : rappelons que l'on a supposé $\tau_0$ unitaire, et donc $\pi$ est 
 de longueur finie et unitaire, donc semi-simple,  et bien sûr, le module standard en position de Langlands négative admet un 
 unique sous-module   irréductible.
 
 Pour cela, nous allons avoir besoin de quelques considérations préliminaires sur les représentations induites des groupes généraux linéaires.
Utilisons les  notations usuelles pour les induites depuis les sous-groupes paraboliques
standard dans les groupes   généraux linéaires et classiques. Comme précédemment, notons 
$\delta_{t,s}$ la représentation de $\GL_2(\bbR)$ essentiellement de carré intégrable
 de caractère infinitésimal $\left( \frac{-t+s}{2},  \frac{t+s}{2}\right)$, où 
$t\in \bbZ_{>0}$ et $s\in \bbC$ (jusque là, nous n'avions introduit que les séries discrètes, c'est-à-dire
les $\delta_{t,s}$ avec $s\in  i \bbR)$, et notons $\eta_{\epsilon, s}$ le caractère de $\GL_1(\bbR)$
défini par $x\mapsto \sgn(x)^{\frac{1-\epsilon}{2}}  \vert x\vert^{s}$, $\epsilon\in {\pm 1}$, $s\in \bbC$, de caractère infinitésimal $s$.
Si $\delta$ est l'une de ces représentations $\delta_{t,s}$ ou $\eta_{\epsilon, s}$ de  $\GL_2(\bbR)$ ou  $\GL_1(\bbR)$, 
notons $\bbZ(\delta)$ l'ensemble $ \frac{t+s}{2}+\bbZ$ dans le premier cas, et $s+\bbZ$ dans le second.
D'après un résultat de B. Speh (\cite{Speh}), si $\delta_1$ et $\delta_2$ sont deux représentations de cette forme, alors 
la représentation induite $\delta_1\times \delta_2$ est irréductible si $\bbZ(\delta_1)\neq \bbZ(\delta_2)$.
En particulier, comme dans le groupe de Grothendieck on a  $\delta_1\times \delta_2=\delta_2\times \delta_1$
en toute généralité, on voit dans ce cas que 
$\delta_1\times \delta_2$ et $\delta_2\times \delta_1$ sont isomorphes.
On a des familles d'opérateurs d'entrelacements 
\[ M(y_1,y_2)  : \delta_1 \vert . \vert  ^{y_1}  \times   \delta_2  \vert . \vert  ^{y_2}  \longrightarrow   \delta_2  \vert . \vert  ^{y_2} \times  \delta_1 \vert . \vert  ^{y_1}\]
et 
\[ N(y_1,y_2)  :   \delta_2  \vert . \vert  ^{y_2} \times  \delta_1 \vert . \vert  ^{y_1} \longrightarrow  \delta_1 \vert . \vert  ^{y_1}  \times   \delta_2  \vert . \vert  ^{y_2}\]
méromorphes en $(y_1,y_2)\in \bbC^2$. Les compositions $M(y_1,y_2)\times N(y_1,y_2)$ et 
$N(y_1,y_2)\times N(y_1,y_2)$ sont des opérateurs scalaires donnés par une même fonction méromorphe 
$\eta(y_1,y_2)$ à valeurs complexes.
Cette fonction n'a pas de pôle en $(0,0)$ (à cause de la condition sur $\bbZ(\delta_1)$ et  $\bbZ(\delta_2)$)    et l'un des opérateurs $M(y_1,y_2)$
ou $N(y_1,y_2)$ est holomorphe en $(0,0)$ (celui pour lequel le terme de gauche de la flèche est en position de Langlands positive, et le terme de droite
en position négative). Les opérateurs d'entrelacement $M(0,0)$ et $N(0,0)$ sont donc définis et réalisent
l'isomorphisme entre $\delta_1\times \delta_2$ et $\delta_2\times \delta_1$.
 
 On tire de ceci le résultat suivant
 \begin{prop} \label{Spehplus} Soient $\tau_1$ et $\tau_2$ des représentations irréductibles de $\GL_{n_1}(\bbR)$ et  $\GL_{n_2}(\bbR)$
 respectivement. Supposons que pour un $x_1\in \bbC$, le caractère infinitésimal de $\tau_1$ soit formé de nombres
 tous dans $x+\bbZ$, et supposons qu'aucune composante du caractère infinitésimal de $\tau_2$ ne soit dans $x+\bbZ$.
 Alors les représentations induites $\tau_1\times \tau_2$ et $\tau_2\times \tau_1$ de $\GL_{n_1+n_2}(\bbR)$ sont isomorphes.
 \end{prop}
 
\dem  On réalise $\tau_1$ et $\tau_2$ comme sous-module de représentations standards en position de Langlands négative :  
 \[\tau_1\hookrightarrow \delta_{1,1}\times \ldots \times \delta_{1,r_1}, \quad \tau_2 \hookrightarrow \delta_{2,1}\times \ldots \times \delta_{2,r_2}.\]
 On a une famille méromorphe d'opérateurs d'entrelacement $\caM(y_1,y_2)$ : 
 \[ \delta_{1,1}\vert . \vert  ^{y_1}\times \ldots \times \delta_{1,r_1}\vert . \vert  ^{y_1}  \times  \delta_{2,1}\vert . \vert  ^{y_2}\times \ldots \times
  \delta_{2,r_2} \vert . \vert  ^{y_2} \longrightarrow   \delta_{2,1}\vert . \vert  ^{y_2}\times \ldots \times
  \delta_{2,r_2} \vert . \vert  ^{y_2} \times  \delta_{1,1}\vert . \vert  ^{y_1} \times \ldots \times \delta_{1,r_1}\vert . \vert  ^{y_1}  \]
 qui se factorise en un produit de composition d'opérateurs d'entrelacement élémentaires de la forme considérée avant l'énoncé de la proposition et qui sont tous holomorphes 
 bijectifs en $(0,0)$. L'opérateur $\caM(y_1,y_2)$ est donc holomorphe bijectif en $(0,0)$. 
\qed

 \ 
 
 Revenons maintenant à notre but principal dans cette cinquième étape.
 
 Pour tout indice $i$, la représentation de  $\mathbf{Speh}\, (\delta_{t_i,s_i},a_i)$ est réalisée comme unique sous-module
 de la représentation standard   \[  S_{t_i,s_i,a_i} = \delta_{t_i,s_i}\  \vert . \vert^{-\frac{a_i-1}{2}} \times   \delta_{t_i,s_i}  \vert . \vert^{-\frac{a_i-3}{2}} \times 
\ldots \times  \delta_{t_i,s_i} \vert .  \vert^{ \frac{a_i-3}{2}}  \times  \delta_{t_i,s_i} \vert . \vert^{\frac{a_i-1}{2}} . \]
 
 Réalisons aussi la représentation $\tau_0$ comme sous-module  de Langlands 
d'une représentation standard $\tilde \tau$ de  $G_0$ en position négative. 
 Ecrivons  $\tilde \tau=\tilde \tau_{--} \times \tilde \tau_{temp}$, où 
 $\tilde \tau_{temp}$ est tempérée irréductible, et $\tau_{--}$ est en position de Langlands strictement négative
 (ici $\tilde \tau_{--}$ est donc une représentation d'un produit de $\GL_1(\bbR)$ et $\GL_2(\bbR)$, et  $\tilde \tau_{temp}$ une représentation d'un produit de
 $\GL_1(\bbR)$ et $\GL_2(\bbR)$ et d'un groupe classique, le  produit de ces deux facteurs   formant un sous-groupe de Levi standard de $G_0$).
 On obtient donc un plongement 
 \begin{equation} \label{Plong} \pi \hookrightarrow  \left(  \times_i   S_{t_i,s_i,a_i} \right) \times   \tilde \tau_{--} \times \tilde \tau_{temp}. \end{equation}

 Formons maintenant une  représentation standard  en position négative de la manière suivante. La représentation  $\times_i   S_{t_i,s_i,a_i}\times \tilde \tau_{--}$
 s'écrit comme un produit de représentations de la forme  $\delta\vert .\vert ^{x}$  où $\delta$ 
 est une série discrète de  $\GL_1(\bbR)$ ou  $ \GL_2(\bbR)$ et 
 $x$ est un demi-entier.
 Remplaçons dans ce produit les termes comme ci-dessus avec $x>0$ par $\delta \vert .\vert ^{-x}$
  et réordonnons les facteurs pour les mettre dans l'ordre des $x$ croissants.
 Ecrivons le produit obtenu comme étant $\Delta_{--}\times \Delta_{temp}$ où $\Delta_{--}$ est le produit des  facteurs $\delta\vert .\vert ^{x}$
 avec $x<0$ et $\Delta_{temp}$ celui avec $x=0$. Notons $\caN$ l'opérateur d'entrelacement standard pour le groupe $G$  qui envoie 
 $   \left(  \times_i   S_{t_i,s_i,a_i} \right) \times   \tau_{--} \times \tilde \tau_{temp}$ dans
  $   \Delta_{--}\times \Delta_{temp}  \times    \tilde \tau_{temp}$. L'opérateur 
  d'entrelacement $\caN$ 
 se factorise en opérateurs élémentaires, l'effet d'un opérateur élémentaire étant de remplacer un facteur $\delta\vert .\vert^x$ 
 avec $x>0$ par     $\delta\vert .\vert^{-x}$ ou bien un produit de la forme $\delta_1\vert .\vert^{x_1}\times \delta_2 \vert .\vert^{x_2}$
 avec $0\geq x_1>x_2$ par  $\delta_2 \vert .\vert^{x_2}\times \delta_1\vert .\vert^{x_1}$ et ceux-ci sont bien définis (holomorphes)
 dans le domaine où on les considère. 
    Notons  encore $\caN$ la composition de $\caN$ avec (\ref{Plong})
  \begin{equation} \label{Plong3} \caN : \, \pi \hookrightarrow  \left(  \times_i   S_{t_i,s_i,a_i} \right) \times \tilde \tau_{--}  \times \tilde \tau_{temp}
  \longrightarrow   \Delta_{--}\times \Delta_{temp}  \times   \tilde \tau_{temp} .\end{equation}
La représentation  $\Delta_{temp}  \times   \tilde \tau_{temp}$ est une représentation tempérée  d'un groupe classique, induite d'une tempérée
irréductible. La théorie du $R$-groupe et l'hypothèse sur les parités de $\Delta_{temp}$ (mauvaise)  et de $\tilde \tau_{temp}$ (bonne)
nous dit que cette représentation est irréductible.
D'autre part, $\Delta_{--}$ est en position de Langlands strictement négative. Le terme de droite est donc une représentation 
standard en position de Langlands  négative, qui admet un unique sous-module irréductible.

Ainsi, on a bien construit un opérateur d'entrelacement  de $\pi$ vers un module standard en position de Langlands négative, et il reste à montrer
son  injectivité.  
Faisons tout d'abord quelques observations sur  $\Delta_{--}$  et $\Delta_{temp}$. 
La première est formée à partir de facteurs   $\delta \vert . \vert ^x$ venant soit des $S_{t_i,s_i,a_i}$, soit de $\tilde \tau_{--}$.
Mais un facteur  $\delta_1 \vert . \vert ^{x_1}$ venant d'un  $S_{t_i,s_i,a_i}$ et 
un facteur $\delta_1 \vert . \vert ^{x_1}$ venant de  $\tilde \tau_{--}$ commutent, car leur produit est irréductible 
d'après le  résultat de B. Speh  et les hypothèses sur les parités.
On a donc en fait  $\Delta_{--}=\Delta_{--}^{Speh}\times \tilde \tau_{--}=   \tilde \tau_{--}\times\Delta_{--}^{Speh} $ où $\Delta_{--}^{Speh}$ est obtenue comme ci-dessus
en changeant des exposants en leurs opposés et en remettant  le tout dans l'ordre, mais seulement pour les facteurs
provenant des représentations de Speh. Le terme $\Delta_{temp}$ est lui un produit de facteurs provenant des représentations de Speh.
On peut donc noter $\Delta_{temp}=\Delta_{temp}^{Speh}$ pour insister sur ce fait. D'autre part, il commute avec $\tilde \tau_{--}$ en vertu du résultat
de Speh invoqué ci-dessus.
Ainsi (\ref{Plong3}) peut aussi s'écrire 
 \begin{equation} \label{Plong6} \caN : \, \pi 
   \longrightarrow   \Delta^{Speh}_{--}\times \Delta^{Speh}_{temp}  \times  \tilde \tau_{--}\times \tilde \tau_{temp} .\end{equation}

Montrons maintenant que pour l'injectivité de $\caN$, on se ramène au  cas où $\tau_0$ est une représentation d'un groupe compact.
En effet, supposons que $\tau_0$ soit sous-module d'une  série principale $(\times_j \gamma_j)\times \tau'_0$ où les $\gamma_j$ 
sont des caractères de $\GL_1(\bbR)$ et $\tau'_0$ est une représentation d'un groupe compact  $G'_0$ de même type que $G$
  ($G'_0$ est la partie compacte du facteur de Levi d'un parabolique minimal de $G_0$). 
On a donc 
\begin{equation} \label{Plong4}\pi=( \times_i  \mathbf{Speh}(\delta_{t_i,s_i},a_i)) \times \tau_0 \hookrightarrow \left( \times_i  \mathbf{Speh}(\delta_{t_i,s_i},a_i)\right)
 \times  (\times_j \gamma_j)\times \tau'_0. \end{equation}
 Grâce à la proposition ci-dessus, $   \left( \times_i  \mathbf{Speh}(\delta_{t_i,s_i},a_i)\right)\times  (\times_j \gamma_j)=
 \times  (\times_j \gamma_j) \times  \left( \times_i  \mathbf{Speh}(\delta_{t_i,s_i},a_i)\right)$
  et l'on   peut réécrire le terme de droite en permutant les facteurs. On obtient donc un plongement 
\begin{equation} \label{Plong7}\pi=( \times_i  \mathbf{Speh}(\delta_{t_i,s_i},a_i)) \times \tau_0 \hookrightarrow
 (\times_j \gamma_j) \times  \left( \times_i  \mathbf{Speh}(\delta_{t_i,s_i},a_i)\right)\times \tau'_0=  (\times_j \gamma_j) \times \pi'. \end{equation}

Admettons le résultat pour $\tau'_0$, à savoir que 
l'opérateur d'entrelacement 
\[\caN' : \,  \pi'= \left( \times_i  \mathbf{Speh}(\delta_{t_i,s_i},a_i)\right)\times \tau'_0  \hookrightarrow    \Delta'_{--}\times \Delta'_{temp}  \times   \tilde \tau'_{temp}  \] 
construit comme ci-dessus en partant de $\pi'$ plutôt que de $\pi$, et avec les notations évidentes, est injectif.
Remarquons que comme $\tau'_0$ est une représentation irréductible d'un groupe compact, on a avec ces notations
$\tilde \tau'=\tilde \tau'_{temp}=\tau'_0$ et $\tilde \tau'_{--}$ est triviale et en particulier $\Delta'_{--}=\Delta^{Speh}_{--}$
D'autre part  $\Delta'_{temp}=\Delta_{temp}^{Speh}=\Delta_{temp}$. On a donc un opérateur injectif 
\[\caN' : \,  \pi'= \left( \times_i  \mathbf{Speh}(\delta_{t_i,s_i},a_i)\right)\times \tau'_0  \hookrightarrow    \Delta^{Speh}_{--}\times \Delta^{Speh}_{temp}  \times    \tau'_{0}.  \] 

Par exactitude du foncteur d'induction parabolique, on en déduit un opérateur  injectif, 
\[\caN' : \,  (\times_j \gamma_j)\times  \pi' \hookrightarrow 
 (\times_j \gamma_j)\times  \Delta^{Speh}_{--}\times \Delta^{Speh}_{temp}  \times    \tau'_{0} . \]
On utilise à nouveau la proposition \ref{Spehplus}
 pour écrire le terme de droite sous la forme  $\Delta^{Speh}_{--}\times \Delta^{Speh}_{temp}  \times   (\times_j \gamma_j)\times   \tau'_{0}$, et 
on compose avec le plongement (\ref{Plong7}) pour obtenir un morphisme injectif  que l'on note encore $\caN'$: 
 \[\caN' : \, \pi \hookrightarrow 
  \Delta^{Speh}_{--}\times \Delta^{Speh}_{temp}  \times   (\times_j \gamma_j)\times   \tau'_{0} . \]
 
 L'injection $\tau_0\hookrightarrow   (\times_j \gamma_j)\times \tau'_0$ est obtenue en composant  l'injection 
 $\tau_0\hookrightarrow \tilde \tau_{--}\times \tilde \tau_{temp} $ et  un morphisme $\tilde \tau_{--}\times \tilde \tau_{temp} \rightarrow  (\times_j \gamma_j)\times
 \tau'_0$.
 Ainsi l'on voit que $\caN'$ se factorise par $\caN$. Ceci établit le fait que  l'injectivité de $\caN'$ implique celle de $\caN$.

Pour les groupes symplectiques,  qui sont déployés,   $\tau'_0$ est la représentation triviale du groupe trivial.
 L'injectivité de $\caN'$ provient alors de l'irréductibilité de $\pi'$ établie à la quatrième étape, et du fait que $\caN'$ est non nul.

 Les groupes orthogonaux demandent encore un peu de travail
 pour conclure comme ci-dessus car $G'_0$ peut être de rang trop petit pour pouvoir appliquer l'irréductibilité 
   démontrée à la quatrième étape. On utilise les   foncteurs de translation pour se ramener 
 au cas où  le caractère infinitésimal de $\tau'_0$ est celui de la représentation triviale, c'est-à-dire que $\tau'_0$ est la représentation
  triviale  $\Triv_{G'_0}$ de $G'_0$.

  Soit $\sigma$ un facteur  irréductible de $\pi'=\left( \times_i  \mathbf{Speh}(\delta_{t_i,s_i},a_i)\right)\times \Triv_{G'_0} $
   (rappelons que $\pi'$ est unitaire et donc semi-simple).
   Notons $n_0$ le plus grand entier entrant dans le caractère infinitésimal de $\Triv_{G'_0}$.
  Soit $T$ un entier suffisamment grand. Notons  $G''_0$ le  groupe  de même type que $G'_0$ et de rang $T+\mathrm{rg}(G'_0)$
  tel que $\GL_1(\bbR)^T\times G''_0$ soit un sous-groupe de Levi de $G'_0$,  et soit 
  $\Triv_{G''_0}$ la représentation triviale de ce groupe.
  On suppose donc $T$ assez  grand pour que les hypothèses de la quatrième étape soient 
  vérifiées pour    $\pi''=\left( \times_i  \mathbf{Speh}(\delta_{t_i,s_i},a_i)\right)\times \Triv_{G''_0} $ qui est donc irréductible.
  La représentation induite 
  \begin{equation}\label{I3} \vert .\vert^{n_0+T}\times \vert .\vert^{n_0+T-1}\times \cdots \times  \vert .\vert^{n_0+1}\times \sigma   \end{equation}
  possède un seul quotient irréductible car elle est quotient de la représentation 
  \[  \vert .\vert^{n_0+T}\times \vert .\vert^{n_0+T-1}\times \cdots \times  \vert .\vert^{n_0+1}\times S(\sigma)   \]
où $S(\sigma)$ est une représentation standard en position de Langlands positive et dont $\sigma$ est le quotient de Langlands.
Grâce à la proposition \ref{Spehplus}, on peut mettre (\ref{I3}) en position de Langlands  positive, et elle admet donc un unique quotient irréductible.  
 Ce quotient irréductible est isomorphe à l'image de l'opérateur d'entrelacement standard
 \[  \vert .\vert^{n_0+T}\times \vert .\vert^{n_0+T-1}\times \cdots \times  \vert .\vert^{n_0+1}\times \sigma 
 \longrightarrow  \vert .\vert^{-n_0-T}\times \vert .\vert^{-n_0-T+1}\times \cdots \times  \vert .\vert^{-n_0-1}\times \sigma  \]
 Donc les quotients irréductibles de 
 \begin{equation}\label{I4} \vert .\vert^{n_0+T}\times \vert .\vert^{n_0+T-1}\times \cdots \times  \vert .\vert^{n_0+1}\times \pi'   \end{equation}
sont isomorphes aux sous-modules    irréductibles de   l'image de l'opérateur d'entrelacement standard
 \begin{equation}\label{E2}  \vert .\vert^{n_0+T}\times \vert .\vert^{n_0+T-1}\times \cdots \times  \vert .\vert^{n_0+1}\times \pi'
 \longrightarrow  \vert .\vert^{-n_0-T}\times \vert .\vert^{-n_0-T-1}\times \cdots \times  \vert .\vert^{-n_0-1}\times \pi'  \end{equation}
et il y a bijection entre ces sous-modules  irréductibles de l'image (et cette image est semi-simple)
et les sous-modules irréductibles  de $\pi'$.


 L'opérateur d'entrelacement (\ref{E2})  se réécrit en utilisant la proposition \ref{Spehplus} ci-dessus 
 \begin{equation}\label{E4} \left( \times_i  \mathbf{Speh}(\delta_{t_i,s_i},a_i)\right)
\times  \vert .\vert^{n_0+T}\times \vert .\vert^{n_0+T-1}\times \cdots \times  \vert .\vert^{n_0+1} 
\times \Triv_{G'_0} \longrightarrow\qquad \qquad  \end{equation}
\[\qquad \qquad \qquad  \left( \times_i  \mathbf{Speh}(\delta_{t_i,s_i},a_i)\right)
\times \vert .\vert^{-n_0-T}\times \vert .\vert^{-n_0-T+1}\times \cdots \times  \vert .\vert^{-n_0-1}\times  \Triv_{G'_0}\]
Or $\vert .\vert^{n_0+T}\times \vert .\vert^{n_0+T-1}\times \cdots \times  \vert .\vert^{n_0+1} 
\times \Triv_{G'_0}$, (resp. $\vert .\vert^{-n_0-T}\times \vert .\vert^{-n_0-T+1}\times \cdots \times  \vert .\vert^{-n_0-1}\times  \Triv_{G'_0}$)
est la représentation standard en position de Langlands positive (resp. négative) dont 
$\Triv_{G''_0}$ est l'unique quotient irréductible (resp. sous-module).
 Ainsi $\pi''=\left( \times_i  \mathbf{Speh}(\delta_{t_i,s_i},a_i)\right)\times \Triv_{G''_0}$
apparaît comme image de (\ref{E4}),  et cette image est irréductible.
Mais (\ref{E4}) s'écrit comme la composition de l'isomorphisme 
\[ \left( \times_i  \mathbf{Speh}(\delta_{t_i,s_i},a_i)\right)
\times  \vert .\vert^{n_0+T}\times \vert .\vert^{n_0+T-1}\times \cdots \times  \vert .\vert^{n_0+1} 
\times \Triv_{G'_0}  \qquad \qquad \qquad \] 
\[  \qquad \qquad \qquad \simeq   \vert .\vert^{n_0+T}\times \vert .\vert^{n_0+T-1}\times \cdots \times  \vert .\vert^{n_0+1} 
\times \left( \times_i  \mathbf{Speh}(\delta_{t_i,s_i},a_i)\right) \times \Triv_{G'_0}
\]
\[ \qquad \qquad \qquad  \qquad \qquad \qquad =   \vert .\vert^{n_0+T}\times \vert .\vert^{n_0+T-1}\times \cdots \times  \vert .\vert^{n_0+1} 
\times \pi',\]
 de (\ref{E2}), et de l'isomorphisme 
\[ \vert .\vert^{n_0+T}\times \vert .\vert^{-n_0-T+1}\times \cdots \times  \vert .\vert^{-n_0-1}\times \pi' \qquad \qquad \qquad  \qquad \qquad \qquad\]
\[\qquad \qquad \qquad =
\vert .\vert^{-n_0-T}\times \vert .\vert^{-n_0-T+1}\times \cdots \times  \vert .\vert^{-n_0-1}\times   \left( \times_i  \mathbf{Speh}(\delta_{t_i,s_i},a_i)\right)\times 
\Triv_{G'_0} \]
\[\qquad \qquad \qquad \simeq 
\left( \times_i  \mathbf{Speh}(\delta_{t_i,s_i},a_i)\right)
\times \vert .\vert^{-n_0-T}\times \vert .\vert^{-n_0-T+1}\times \cdots \times  \vert .\vert^{-n_0-1}\times  \Triv_{G'_0}\]

Ainsi l'image de (\ref{E4}) est d'une part irréductible et d'autre part a autant de composantes irréductibles que $\pi'$  d'après la remarque faite après  (\ref{E2}).
 Ceci montre que $\pi'$ est irréductible.
\

 Ceci termine la démonstration de la proposition  \ref{redCla2}.\qed 
 
\bigskip
\bibliographystyle{smfalpha}
\bibliography{MR7}

  \end{document}